\documentclass[11pt,a4paper]{article} 
\usepackage{url}

\usepackage{tikz}
\tikzset{every picture/.style={line width=0.75pt}} 

 \usepackage[usenames,dvipsnames]{pstricks}
 \usepackage{epsfig}
 \usepackage{pst-grad} 
 \usepackage{pst-plot} 
 \usepackage[space]{grffile} 
 \usepackage{etoolbox} 
 \makeatletter 
 \patchcmd\Gread@eps{\@inputcheck#1 }{\@inputcheck"#1"\relax}{}{}
 \makeatother

\usepackage[T1]{fontenc} 

\usepackage{amsmath,mathtools}
\usepackage{amssymb}
\usepackage{amsthm}
\usepackage[latin1]{inputenc}
     
\usepackage{graphicx}

\usepackage[usenames,dvipsnames]{pstricks}
\usepackage{epsfig}
\usepackage{pst-grad} 
\usepackage{pst-plot} 

\newcommand\restr[2]{{
  \left.\kern-\nulldelimiterspace 
  #1 
  \vphantom{\Big|} 
  \right|_{#2} 
  }}

\usepackage[misc,geometry]{ifsym}

\usepackage{ifthen}
\usepackage{color}

\newcommand{\intav}[1]{\mathchoice {\mathop{\vrule width 6pt height 3 pt depth  -2.5pt
\kern -8pt \intop}\nolimits_{\kern -6pt#1}} {\mathop{\vrule width
5pt height 3  pt depth -2.6pt \kern -6pt \intop}\nolimits_{#1}}
{\mathop{\vrule width 5pt height 3 pt depth -2.6pt \kern -6pt
\intop}\nolimits_{#1}} {\mathop{\vrule width 5pt height 3 pt depth
-2.6pt \kern -6pt \intop}\nolimits_{#1}}}

\def\polhk#1{\setbox0=\hbox{#1}{\ooalign{\hidewidth\lower1.5ex\hbox{`}\hidewidth\crcr\unhbox0}}}

\def\XXint#1#2#3{{\setbox0=\hbox{$#1{#2#3}{\int}$ }
\vcenter{\hbox{$#2#3$ }}\kern-.6\wd0}}

\usepackage{setspace}

\newtheorem{theorem}{Theorem}
\newtheorem{example}{Example}
\newtheorem{definition}{Definition}
\newtheorem{lemma}{Lemma}
\newtheorem{corollary}{Corollary}
\newtheorem{proposition}{Proposition}
\newtheorem{remark}{Remark}
\newtheorem{assumption}{A}

\definecolor{r2}{RGB}{0, 0, 0}
\definecolor{r1}{RGB}{0, 0, 0}
\definecolor{wethepeople}{RGB}{211, 211, 211}
\usepackage[table]{xcolor}

\makeatletter
\patchcmd{\env@cases}{1.2}{1}{}{}
\makeatother

\begin{document}

\title{Numerical methods for fully nonlinear degenerate diffusions}

\author{Edgard A. Pimentel and Erc\'ilia Sousa}
	
\date{\today} 

\maketitle

\begin{abstract}
We propose finite difference methods for degenerate fully nonlinear elliptic equations and prove the convergence of the schemes. Our focus is on the pure equation and a related free boundary problem of transmission type. The cornerstone of our argument is a regularisation procedure. It decouples the degeneracy term from the elliptic operator driving the diffusion process. In the free boundary setting, a two-parameters regularisation entails new, genuine difficulties. To bypass them, we resort to the intrinsic properties of the regularised problem and a fixed-point type of argument. We present numerical experiments supporting our theoretical results. Our methods are flexible, and our approach can be extended to a broader class of non-variational problems.
\end{abstract}

\smallskip


\smallskip 


\vspace{.1in}

\section{Introduction}\label{sec_intro}

We propose numerical methods for fully nonlinear degenerate equations of the form
\begin{equation}\label{eq_main1}
	\begin{cases}
	\left|Du\right|^\theta F(D^2u)=f&\hspace{.2in}\mbox{in}\hspace{.1in}\Omega\\
	u=g&\hspace{.2in}\mbox{on}\hspace{.1in}\partial\Omega,
	\end{cases}
\end{equation}
where $F:S(d)\to\mathbb{R}$ is a fully nonlinear elliptic operator, $f\in C(\Omega)\cap L^\infty(\Omega)$, and $g\in C(\partial\Omega)$. The exponent $\theta>0$ is a degeneracy rate, governing how the vanishing of the gradient affects the ellipticity of the problem.

After designing a finite difference  scheme for \eqref{eq_main1}, we embed this equation into a free boundary problem. Namely, we consider the free transmission model
\begin{equation}\label{eq_main2}
	\begin{cases}
		|Du|^{\theta_1}F(D^2u)=f&\hspace{.2in}\mbox{in}\hspace{.1in}\Omega\cap\left\lbrace u>0\right\rbrace\\
		|Du|^{\theta_2}F(D^2u)=f&\hspace{.2in}\mbox{in}\hspace{.1in}\Omega\cap\left\lbrace u<0\right\rbrace\\
		u=g&\hspace{.2in}\mbox{on}\hspace{.1in}\partial\Omega,
	\end{cases}
\end{equation}
where {\color{r2}$0\leq\theta_1<\theta_2$}. In both settings, the unknown $u$  is understood in the sense of viscosity solutions \cite{CraLio,CraEvaLio,CraIshLio,CafCab}.

The motivation to enhance the understanding of \eqref{eq_main1} and \eqref{eq_main2} is twofold. First, this class of equations generalise usual structural conditions associated with fully nonlinear elliptic problems. In addition, the uniqueness of solutions to \eqref{eq_main1} and \eqref{eq_main2} is unknown in general, and numerical approximations provide insights on selection criteria (c.f. vanishing viscosity methods). Also, \eqref{eq_main2} cannot be written as a single equation in $\Omega$ (unlike the obstacle problem), and an alternative approach is required. 

When it comes to applications, several models can be cast as \eqref{eq_main1} or in a similar fashion. Indeed, the PDE in \eqref{eq_main1} accounts for a fully nonlinear counterpart of the $p$-Laplacian, as the game-theoretical approach to that operator leads to
\begin{equation}\label{eq_dl}
	\Delta_p u=|Du|^{p-2}\left(\Delta u+(p-2)\left\langle D^2u\frac{Du}{|Du|},\frac{Du}{|Du|}\right\rangle\right);
\end{equation}
see \cite{PerSchSheWil}. Although the latter and \eqref{eq_main1} involve distinct structural conditions, they share important similarities, most notably in terms of the product structure and monotonicity properties. One can also formulate {\color{r1} the diffusion part of a }non-local porous medium equations in the form \eqref{eq_main1}, at least in dimension $d=1$. {\color{r1}Indeed, for $s\in(0,1)$, consider the scalar parabolic model
\begin{equation*}\label{eq_parpme}
	\partial_tu=\partial_x\left(u^{m-1}\partial_x(-\Delta)^{-s}u\right),
\end{equation*}
with $m\in(1,\infty)$, $x\in\mathbb{R}$, and $t>0$. In \cite{StadelVaz}, the authors use the cumulative density $v=v(x,t)\geq 0$, that satisfies $v_x=u$, and 
\[
	\partial_tv=-|v_x|^{m-1}\left(-\Delta\right)^\alpha v,
\]
in the viscosity sense, where $\alpha=1-s$. For non-variational formulations of the porous medium equation, degenerating with a solution-dependent rate, we refer the reader to \cite{BraVaz} and \cite{ChaSan}.
} 

Finally, non-local quasilinear equation may be reduced (under suitable conditions) to a non-local variant of the model in \eqref{eq_main1}; see \cite{ChaJak}. We mention that \eqref{eq_main1} has been largely studied in the literature; see, for instance, \cite{BirDem2004,DavFelQua}. Regularity results for the solutions to \eqref{eq_main1} are the subject of \cite{ImbSil2013,BirDem2014}. The transmission problem in \eqref{eq_main2} was introduced in \cite{HuaPimRamSwi}.

Numerical methods for nonlinear partial differential equations (PDE) of Hamilton-Jacobi type have appeared in \cite{BarFal,CapFal,Sou}, to name just a few. In \cite{BarSou}, the authors examine a strategy for the approximation of viscosity solutions to fully nonlinear equations. They also devise a set of conditions under which a family of numerical approximations converge to the (unique) viscosity solution of the underlying PDE. The authors prove that solutions to a monotone, consistent and stable scheme converge to the unique viscosity solution of the approximated problem. In \cite{Obe}, the author verifies that degenerate ellipticity, properness and Lipschitz continuity of the numerical method ensure monotonicity and stability and lead to convergence. See also \cite{BonZid,Fenjen,DebJak,ObeSal2015,FroObeSal,CagGomTra,Bonnet}.

{\color{r1}We note there has been also recent developments on finite difference schemes concerning degenerate diffusions driven by the $p$-Laplace operator. In \cite{delLin_2022} the authors propose a finite differences discretisation for the Dirichlet problem governed by \eqref{eq_dl}. Their techniques are cleverly based on a new mean-value formula for the $p$-Laplace operator (of fundamental, independent interest in itself; see for instance \cite{delRos_2026} and the references therein). The parabolic variant of the $p$-Poisson equation is the subject of \cite{delLin_2023}, where the authors develop a finite difference scheme also based on the mean-value property. We do not expect the strategy in those papers would extend to the fully nonlinear setting. On the other hand, we believe the techniques in the present paper could provide a companion approach to the methods in \cite{delLin_2022}. If this is the case, a comparative analysis concerning the computational cost of both schemes would play an important role; we do not pursue this endeavour in the present paper.
}

Solving a nonlinear approximation scheme designed for a stationary PDE benefits from a solution operator in the spirit of the Euler method. In \cite{Obe}, the author resorts to this method and derives a non-linear Courant-Friedrichs-Lewy (CFL) condition ensuring the solution operator is a contraction on a Banach space. Such a condition depends only on the Lipschitz constant of the approximation scheme.  

Problems \eqref{eq_main1} and \eqref{eq_main2} introduce genuine difficulties concerning numerical approximations. The first one regards monotonicity. It is well-known that there exists a monotone approximation for $F(D^2u)$, provided $F$ satisfies a diagonal-dominance condition (see \cite[Theorem 3.4]{BarSou}). On the other hand, an upwind monotone discretisation of $Du$ yields a monotone approximation for $|Du|^\theta$. However, the \emph{product} of monotone operators are not necessarily monotone.

The second main difficulty concerns the (two-phase) free boundary problem \eqref{eq_main2}. Unlike the obstacle problem, \eqref{eq_main2} cannot be immediately written in terms of a single equation in the entire domain. Compare with the \cite[Section 3]{Obe}. Therefore, an approximation scheme in $\Omega$ is not available. In addition, the dependence of the degeneracy rates on the sign of the solution affects the maximum and the comparison principles.

We start our analysis with the pure equation \eqref{eq_main1}. For $0<\varepsilon\ll1$, we propose a regularisation of the form
\begin{equation}\label{eq_reg1}
	\begin{cases}
		\left(\varepsilon+|Du^\varepsilon|^2\right)^\frac{\theta}{2}\left(\varepsilon u^\varepsilon+F(D^2u^\varepsilon)\right)=f&\hspace{.2in}\mbox{in}\hspace{.1in}\Omega\\
		u^\varepsilon=g&\hspace{.2in}\mbox{on}\hspace{.1in}\partial\Omega.
	\end{cases}
\end{equation}

Arguing along the same lines as in \cite{HuaPimRamSwi}, we prove the existence of a (unique) viscosity solution to \eqref{eq_reg1}. We also verify that $u^\varepsilon\to u$, locally uniformly, where $u\in C(\Omega)$ is a viscosity solution to \eqref{eq_main1}. This is summarised in our first result.

\begin{theorem}[Existence of solutions]\label{thm_existence}
Let $\Omega\subset\mathbb{R}^d$ be a bounded domain satisfying a uniform exterior sphere condition. Suppose Assumptions \ref{assump_A40} and \ref{assump_data}, to be detailed further, hold. Then there exists a unique viscosity solution to \eqref{eq_reg1}, denoted with $u^\varepsilon\in C(\Omega)$. Moreover, $u^\varepsilon\to u$ locally uniformly in $\Omega$, where $u\in C(\Omega)$ is a viscosity solution to \eqref{eq_main1}.
\end{theorem}

Thus, our goal is to produce a numerical approximation to \eqref{eq_reg1}. To that end, we fix $0<h\ll1$ and consider a discrete approximation of $\Omega$, which we denote with $\Omega_h$. Then we propose an approximation of the form
\begin{equation}\label{eq_method1}
	\begin{cases}
		\varepsilon u_h^\varepsilon(x)+F_h(D^2_hu^\varepsilon_h(x))-\frac{f(x)}{\left(\varepsilon + |D_hu_h^\varepsilon(x)|^2\right)^\frac{\theta}{2}}=0&\hspace{.2in}\mbox{in}\hspace{.1in}\Omega_h\\
		u^\varepsilon_h(x)=g(x)&\hspace{.2in}\mbox{on}\hspace{.1in}\partial\Omega_h.
	\end{cases}
\end{equation}

We prove the method in \eqref{eq_method1} is monotone, consistent with \eqref{eq_reg1}, and stable. Therefore, an off-the-shelf application of the Barles-Souganidis theory implies that $u_h^\varepsilon\to u^\varepsilon$ locally uniformly \cite{BarSou}. Finally, by sending $\varepsilon\to 0$, one concludes that $u_h^\varepsilon$ provides an approximation to the solutions to \eqref{eq_main1}. 

We emphasise the regularisation in \eqref{eq_reg1} decouples the product in \eqref{eq_main1}, unlocking monotonicity. In addition, the proof of the existence of viscosity solutions to \eqref{eq_reg1} relies on global barriers. We show that a discrete version of those barriers amounts to sub and supersolutions to the numerical method in \eqref{eq_method1}. Such discrete barriers are independent of the mesh size $h$. Together with the degenerate ellipticity of the method, it provides us with the stability for $u^\varepsilon_h$. To prove the consistency, we rely on standard arguments. Our first main result is the following.

{\color{r1}
\begin{theorem}[Convergence of the numerical method I]\label{thm_convergencepure}
Let $\Omega\subset\mathbb{R}^d$ be a bounded domain satisfying a uniform exterior sphere condition. Suppose Assumptions \ref{assump_A40}-\ref{assump_data}, to be detailed further, hold true. Then the family $(u_h^\varepsilon)_{0<h,\varepsilon<1}$ converges locally uniformly to $u^\varepsilon$ as $h\to 0$, where $u^\varepsilon$ is a viscosity solution to \eqref{eq_reg1}. By taking the limit $\varepsilon\to 0$, we have $u^\varepsilon\to u$, subsequentially if necessary, where $u$ is a viscosity solution of \eqref{eq_main1}.
\end{theorem}
}

Theorem \ref{thm_convergencepure} verifies that the method in \eqref{eq_method1} is monotone and consistent. After proving that global barriers for \eqref{eq_reg1} provide us with global barriers for \eqref{eq_method1}, uniformly in $h$, we prove the method is stable. Hence, $u_h^\varepsilon$ converges locally uniformly to the viscosity solution $u^\varepsilon$ of \eqref{eq_reg1}. 

Consequently, the method \eqref{eq_method1} provides a numerical approximation for \eqref{eq_main1}. It remains to notice that the existence of solutions to \eqref{eq_method1} follows from the Euler method, under an appropriate CFL condition.

Once the numerical approximation for the pure equation is understood, we address the free transmission problem \eqref{eq_main2}. Once again, we resort to a regularisation strategy. However, to accommodate a free boundary problem whose degeneracy law depends on the solutions, we start with an auxiliary function. {\color{r2}Let $v\in C(\overline\Omega)$ and define $g^v_\varepsilon:\overline\Omega\to\mathbb{R}$ as 
\[
	g_\varepsilon^v(x)\coloneqq\max\left(\min\left(\frac{v+\varepsilon}{2\varepsilon},1\right),0\right).
\]
Set $h_\varepsilon^v(x)\coloneqq\left(g_\varepsilon^v\ast\eta_\varepsilon\right)(x)$, where $\eta_\varepsilon$ is a standard mollifier. Finally, define the exponent $\theta_\varepsilon^v:\overline\Omega\to\mathbb{R}$ as
\[
	\theta_\varepsilon^v(x)\coloneqq h_\varepsilon^v(x)\theta_1+(1-h_\varepsilon^v(x))\theta_2.
\]
}
{\color{r2}Notice that $\theta_\varepsilon\in[\theta_1,\theta_2]$. Also, if $t>\varepsilon$ we get $\theta_\varepsilon(t)=\theta_1$; if $t<-\varepsilon$, we get $\theta_\varepsilon(t)=\theta_2$. In the interval $(-\varepsilon,\varepsilon)$, the exponent $\theta_\varepsilon(t)$ is smooth. We consider
\begin{equation}\label{eq_reg2}
	\begin{cases}
		\left(\varepsilon+|Du^{\varepsilon,v}|^2\right)^\frac{\theta^v_\varepsilon(x)}{2}\left(\varepsilon u^{\varepsilon,v}+F(D^2u^{\varepsilon,v})\right)=f&\hspace{.2in}\mbox{in}\hspace{.1in}\Omega\\
		u^{\varepsilon,v}=g&\hspace{.2in}\mbox{on}\hspace{.1in}\partial\Omega.
	\end{cases}
\end{equation}
The strategy in \eqref{eq_reg2} relies on a two-parameter regularisation; see \cite{HuaPimRamSwi}. In that paper, the functional parameter $v\in C(\overline\Omega)$ is dealt with through a fixed-point strategy. Ultimately, the authors obtain the Dirichlet problem
\begin{equation}\label{eq_reg3}
	\begin{cases}
		\left(\varepsilon+|Du^{\varepsilon}|^2\right)^\frac{\theta^{u^\varepsilon}_\varepsilon(x)}{2}\left(\varepsilon u^{\varepsilon}+F(D^2u^{\varepsilon})\right)=f&\hspace{.2in}\mbox{in}\hspace{.1in}\Omega\\
		u^{\varepsilon}=g&\hspace{.2in}\mbox{on}\hspace{.1in}\partial\Omega.
	\end{cases}
\end{equation}
Then, compactness allows them to consider the limit $\varepsilon\to 0$, producing a viscosity solution to \eqref{eq_main2}. Our argument is similar, and transpose the two-parameter regularisation to the numerical setting.
%

We propose a discrete approximation for \eqref{eq_reg2} of the form
\begin{equation}\label{eq_method2}
	\begin{cases}
		\varepsilon u_h^{\varepsilon,v}(x)+F_h(D^2_hu^{\varepsilon,v}_h(x))-\frac{f(x)}{\left(\varepsilon + |D_hu_h^{\varepsilon,v}(x)|^2\right)^\frac{\theta_\varepsilon^v(x)}{2}}=0&\hspace{.2in}\mbox{in}\hspace{.1in}\Omega_h\\
		u^{\varepsilon,v}_h(x)=g(x)&\hspace{.2in}\mbox{on}\hspace{.1in}\partial\Omega_h.
	\end{cases}
\end{equation}

Our analysis start by noticing that $v$ restricted to $\Omega_h$ yields a grid function $v_h$ whose image is a vector in $\mathbb{R}^{N(h)}$. For $0<h\ll1$ and $0<\varepsilon\ll1$ fixed, we prove the existence of $u_h^\varepsilon$ solving 

\begin{equation}\label{eq_method3}
	\begin{cases}
		\varepsilon u_h^{\varepsilon}(x)+F_h(D^2_hu^{\varepsilon}_h(x))-\frac{f(x)}{\left(\varepsilon + |D_hu_h^{\varepsilon}(x)|^2\right)^\frac{\theta_\varepsilon^{u_h^\varepsilon}(x)}{2}}=0&\hspace{.2in}\mbox{in}\hspace{.1in}\Omega_h\\
		u^{\varepsilon}_h(x)=g(x)&\hspace{.2in}\mbox{on}\hspace{.1in}\partial\Omega_h.
	\end{cases}
\end{equation}
The proof follows from a fixed-point argument in finite dimensional spaces. At this point, the analysis focuses on the convergence of the solutions $u_h^\varepsilon$ to \eqref{eq_method3} to a viscosity solution of \eqref{eq_main2}. This is the content of our second main result.

%

\begin{theorem}[Convergence of the numerical method II]\label{thm_conv2}
Suppose Assumptions \ref{assump_dsquare}-\ref{assump_theta}, to be detailed further, are in force. Then the method in \eqref{eq_method3} is monotone, consistent with \eqref{eq_reg3}, and stable. In addition, $u^\varepsilon_h\to u^\varepsilon$ locally uniformly as $h\to 0$, where $u^\varepsilon$ is a viscosity solution to \eqref{eq_reg3}. By taking the limit $\varepsilon\to 0$, we have $u^\varepsilon\to u$, subsequentially if necessary, where $u$ is a viscosity solution of \eqref{eq_main2}.
\end{theorem}

The proof of Theorem \ref{thm_conv2} follows along the same arguments as in the proof of Theorem \ref{thm_convergencepure}. The main difference is the fixed-point analysis accounting for the functional parameter $v$ used in the regularisation of \eqref{eq_main2}.
}

{\color{r1}
\begin{remark}[Joint convergence and the $(h,\varepsilon)$-relation]\label{rem_he} 
An important aspect of Theorems \ref{thm_convergencepure} and \ref{thm_conv2} is the fact that the limits $h\to 0$ and $\varepsilon\to 0$ are taken \emph{separately}. A natural question concerns the possibility of taking joint limits $(h,\varepsilon)\to(0,0)$. One naturally considers
\[
	\sup_{x\in\Omega_h}\left|u_h^\varepsilon(x)-u(x)\right|\leq \sup_{x\in\Omega_h}|u_h^\varepsilon(x)-u^\varepsilon(x)|+\sup_{x\in\Omega_h}|u^\varepsilon(x)-u(x)|.
\]
Suppose there are moduli of continuity $\omega_1$ and $\omega_2$ such that 
\[
	\sup_{x\in\Omega_h}|u_h^\varepsilon(x)-u^\varepsilon(x)|\leq \omega_1(h)\hspace{.2in}\mbox{and}\hspace{.2in}	\sup_{x\in\Omega_h}|u^\varepsilon(x)-u(x)|\leq \omega_2(\varepsilon).
\]
In this case, Theorems \ref{thm_convergencepure} and \ref{thm_conv2} could be phrased as joint-convergence results. We notice that obtaining such moduli is tantamount to convergence rates for the numerical method and the regularised approximation. However, these very important questions of the theory are off the scope of the present manuscript and are not pursued here. We include in Section \ref{sec_numerexp} a numerical experiment relating both parameters and comparing the magnitude of the associated errors.
\end{remark}
}

\begin{remark}[Uniqueness and the selection of solutions]\label{rem_selection} 
For every $0<\varepsilon\ll1$, \eqref{eq_method1} and \eqref{eq_method2} have unique solutions, converging to the (unique) viscosity solution of \eqref{eq_reg1} and \eqref{eq_reg2}, respectively. Although the uniqueness of solutions for \eqref{eq_main1} and \eqref{eq_main2} has not been established, one can use the subsequential limits $\varepsilon\to 0$ to select particular solutions to these equations. It would be interesting to understand whether such (families of) solutions have special properties or satisfy particular conditions (such as transmission conditions). We do not pursue these topics in the present paper.
\end{remark}

{\color{r1}
\begin{remark}[Sign-changing source terms]\label{rem_homog}
In case $f$ changes sign in $\Omega$, and vanishes at some points of the domain, the arguments in this paper must be adjusted. We believe this situation can be tackled by reversing the sign in the discretisation of $|Du|^2$ at points where $f<0$ and considering the uniformly elliptic variant of the model at the points where $f=0$.
\end{remark}
}

The remainder of this paper is organised as follows. Section \ref{sec_prelim} gathers preliminary material, whereas Section \ref{sec_newproofeist} examines the existence of solutions to \eqref{eq_reg1}. The proof of Theorem \ref{thm_convergencepure} is the subject of Section \ref{sec_discrete1}. In Section \ref{sec_pure} we detail the proof of Theorem \ref{thm_conv2}. A final section presents several numerical experiments as an attempt to illustrate our strategy and validate our method.

\section{Main assumptions and preliminary material}\label{sec_prelim}

We denote with $S(d)$ the space of symmetric matrices of order $d$, and notice $S(d)\sim \mathbb{R}^\frac{d(d+1)}{2}$. The norm of $M\in S(d)$ is given by.
\[
	\left\|M\right\|\coloneqq{\color{r2}\sup_{e\in \mathbb{S}^{d-1}}\left|e^TMe\right|.}
\]
For $\xi\in\mathbb{R}^d$, we define $\left\|\xi\right\|\coloneqq \sqrt{\xi\cdot\xi}$. The matrix norm $\left\|N\right\|$ is defined as
\[
	\left\|N\right\|:=\max_{i=1,\ldots,d}\left|e_i\right|,
\]
where $\left\lbrace e_1,e_2,\ldots,e_d \right\rbrace$ are the eigenvalues of the matrix $N$. We also notice that the matrix inequality $M\geq N$ means that $M-N$ is positive semidefinite.

Our first main assumption concerns the uniform ellipticity of the fully nonlinear operator driving \eqref{eq_main1} and \eqref{eq_main2}.

\begin{assumption}[Uniform ellipticity]\label{assump_A40}
Let $0<\lambda\leq\Lambda$ be fixed, though arbitrary. We suppose $F$ is uniformly elliptic, or $(\lambda,\Lambda)$-elliptic. Namely, for every $M,N\in S(d)$ with $N\geq 0$, we have
\[
	\lambda\left\|N\right\|\leq F(M)-F(M+N)\leq \Lambda\left\|N\right\|.
\]
\end{assumption}
 
A well-known consequence of uniform ellipticity is the Lipschitz continuity of the operator $F$ concerning the Hessian entry. Also, a uniformly elliptic operator can be written as an Isaacs-type operator \cite[Remark 1.5]{CabCaf}. Building on that observation, we introduce our second assumption.

\begin{assumption}[Diagonal dominance]\label{assump_dsquare}
We suppose the operator $F:S(d)\to\mathbb{R}$ is of Isaacs-type. That is,
\[
	F(M)=\sup_{\alpha\in\mathcal{A}}\inf_{\beta\in\mathcal{B}}\left(-{\rm Tr}\left(A_{\alpha,\beta}M\right)\right),
\]
where the matrices 
\[
	A_{\alpha,\beta}=\left(a_{i,j}^{\alpha,\beta}\right)_{i,j=1}^d\in S(d)
\]
are $(\lambda,\Lambda)$-elliptic and the sets $\mathcal{A}$ and $\mathcal{B}$ are compact. Also, $A_{\alpha,\beta}$ satisfies a diagonal dominance condition of the form
\[
	a_{i,i}^{\alpha,\beta}\geq {\color{r2}\sum_{\substack{j=1\\j\neq i}}^d\left|a_{i,j}^{\alpha,\beta}\right|},
\]
for every $i=1,\ldots,d$, $\alpha\in\mathcal{A}$ and $\beta\in\mathcal{B}$.
\end{assumption}

The diagonal dominance condition in Assumption \ref{assump_dsquare} ensures that there exists a monotone approximation scheme for the operator $F$; see \cite[Theorem 3.4]{BarSou}. We work under Assumption \ref{assump_dsquare} to simplify matters. {\color{r2} For instances where such a condition has been relaxed, we refer the reader to \cite{BonZid} and \cite{DebJak}, to name just a very few.}

\begin{assumption}[Regularity of the data]\label{assump_data}
We suppose $f\in C(\Omega)\cap L^\infty(\Omega)$ {\color{r1}satisfies $f>0$}. We also suppose $g\in C(\partial\Omega)$.
\end{assumption}

When working under Assumption \ref{assump_data}, one may consider $f\equiv 1$, for simplicity. We also impose conditions on the exponents $\theta_1$ and $\theta_2$.

\begin{assumption}[Degeneracy rates]\label{assump_theta}
We suppose {\color{r2}$0\leq\theta_1<\theta_2$}.
\end{assumption}

We continue with the definition of degenerate ellipticity. 

\begin{definition}[Degenerate elliptic operators]\label{def_dgel}
We say the operator $F:S(d)\times\mathbb{R}^d\times\mathbb{R}\times\Omega\to\mathbb{R}$ is degenerate elliptic if $F(N,p,r,x)\leq F(M,p,r,x)$, for every $p\in\mathbb{R}^d$, $r\in\mathbb{R}$ and $x\in\Omega$, whenever $M\leq N$.
\end{definition}

A $(\lambda,\Lambda)$-elliptic operator is degenerate elliptic. Of particular interest for us is the operator 
\[
	{\color{r2}F_\varepsilon(M,p,r,x)\coloneqq \varepsilon r+F(M)-\frac{f(x)}{\left(\varepsilon+|p|^2\right)^\frac{\theta^r_\varepsilon(x)}{2}}.}
\]
Notice $F_\varepsilon$ is degenerate elliptic provided $F$ is degenerate elliptic. 

Let $G:S(d)\times\mathbb{R}^d\times\mathbb{R}\times\Omega\to\mathbb{R}$ be a $(\lambda,\Lambda)$-elliptic operator and consider the Dirichlet problem
\begin{equation}\label{eq_sundayroast}
	\begin{cases}
		G(D^2u,Du,u,x)=f&\hspace{.2in}\mbox{in}\hspace{.2in}\Omega\\
		u=g&\hspace{.2in}\mbox{in}\hspace{.2in}\partial\Omega,
	\end{cases}
\end{equation}
where $f\in L^p(\Omega)\cap C(\Omega)$ and $g\in C(\partial\Omega)$. We define the (perhaps discontinuous) operator $F:S(d)\times\mathbb{R}^{d}\times \mathbb{R}\times\overline\Omega\to\mathbb{R}$ as
\begin{equation*}\label{eq_woodstock}
	F(M,p,r,x):=
		\begin{cases}
			G(M,p,r,x)-f(x)&\hspace{.2in}\mbox{if}\hspace{.2in}x\in\Omega\\
			r-g(x)&\hspace{.2in}\mbox{if}\hspace{.2in}x\in\partial\Omega.
		\end{cases}
\end{equation*}
The operator $F$ simultaneously accounts for the PDE and the boundary condition. To define a viscosity solution to \eqref{eq_sundayroast}, one requires the notion of semicontinuous envelopes, which is the subject of the next definition.
\begin{definition}[Semicontinuous envelopes]\label{def_envelope}
	Let $F:S(d)\times\mathbb{R}^{d}\times \mathbb{R}\times\overline\Omega\to\mathbb{R}$. We define the upper semicontinuous envelope $F^*:S(d)\times\mathbb{R}^{d}\times \mathbb{R}\times\overline\Omega\to\mathbb{R}$ as
	\[
		F^*(M,p,r,x_0):=\limsup_{x\to x_0}F(M,p,r,x).
	\]
	We also define the lower semicontinuous envelope $F_*:S(d)\times\mathbb{R}^{d}\times \mathbb{R}\times\overline\Omega\to\mathbb{R}$ as
	\[
		F_*(M,p,r,x_0):=\liminf_{x\to x_0}F(M,p,r,x).
	\]
\end{definition}

We proceed with the definition of viscosity solution to \eqref{eq_sundayroast} in terms of the semicontinuous envelopes for $F$.

\begin{definition}[Viscosity solution]\label{def_sandysbar}
An upper semicontinuous function $u:\overline\Omega\to\mathbb{R}$ is a viscosity subsolution to \eqref{eq_sundayroast} if for every $\varphi\in C^2(\overline\Omega)$ such that $u-\varphi$ attains a maximum at $x_0\in\overline\Omega$, we have
\[
	F_*(D^2\varphi(x_0),D\varphi(x_0),u(x_0),x_0)\leq 0.
\]
A lower semicontinuous function $u:\overline\Omega\to\mathbb{R}$ is a viscosity supersolution to \eqref{eq_sundayroast} if for every $\varphi\in C^2(\overline\Omega)$ such that $u-\varphi$ attains a minimum at $x_0\in\overline\Omega$, we have
\[
	F^*(D^2\varphi(x_0),D\varphi(x_0),u(x_0),x_0)\geq 0.
\]
A continuous function that is both a subsolution and a supersolution to \eqref{eq_sundayroast} is a viscosity solution to \eqref{eq_sundayroast}.
\end{definition}

Let $0<h\ll1$ be fixed {\color{r2} and define the $h$-lattice $\Gamma_h\coloneqq h\mathbb{Z}^d$. Define $\Omega_h\coloneqq \Omega\cap\Gamma_h$ and denote with $V$ the stencil
\[
	V\coloneqq\left\lbrace \pm e_i,\,\pm(e_i\pm e_j)\,|\,1\leq i<j\leq d,\,\mbox{$e_i$ and $e_j$ are canonical directions in $\mathbb{R}^d$}\right\rbrace.
\]
Define the set of interior points $\Omega^\circ_h\subset\Omega_h$ as
\[
	\Omega^\circ_h\coloneqq\left\lbrace x\in\Omega_h\,|\,x+hv\in\Omega_h,\,\mbox{for every}\, v\in V\right\rbrace.
\]
Lastly, the boundary of $\Omega_h$ is defined as $\partial\Omega_h\coloneqq\Omega_h\setminus\Omega_h^\circ$.
}
For $0<h\ll1$ fixed, the cardinality of $\Omega_h$ is finite, and denoted with $N_h$. For $0<h\ll1$, we use discrete approximations of the Hessian and the gradient of $u$. {\color{r2} The latter is denoted by $|D_hu_h|^2$ and given by
\[
	\left|D_hu_h(x)\right|^2\coloneqq\frac{1}{h^2}\sum_{i=1}^d\max\left(\left(u_h(x)-u_h(x+he_i)\right),\left(u_h(x)-u_h(x-he_i)\right),0\right)^2.
\]
As concerns the discretisation of the Hessian, we work under Assumption \ref{assump_dsquare}. First, define $\delta_h^{i,i}u_h$ and $\delta_h^{i,j}u_h$ as
\[
	\delta_h^{i,i}u_h(x)\coloneqq\frac{2u_h(x)-u_h(x+e_ih)-u_h(x-e_ih)}{h^2}
\]
and
\[
	\delta_h^{i,j}u_h(x)\coloneqq\frac{2u_h(x)-u_h(x+e_ih+\operatorname*{sgn}(a_{i,j}^{\alpha,\beta})e_jh)-u_h(x-e_ih-\operatorname*{sgn}(a_{i,j}^{\alpha,\beta})e_jh)}{h^2},
\]
for $i,j=1,\ldots,d$, where $A_{\alpha,\beta}$ is the matrix in Assumption \ref{assump_dsquare}. Then we approximate $F(D^2u)$ by $F_h(D^2_hu_h)$, where the latter is defined as
\begin{equation*}
	\begin{split}
		F_h(D^2_hu_h)&\coloneqq \sup_{\alpha\in\mathcal A} \inf_{\beta\in\mathcal B}\sum_{i=1}^d\left(a_{ii}^{\alpha,\beta}-\sum_{j\neq i}|a_{ij}^{\alpha,\beta}|\right)\delta_h^{i,i}u_h(x)\\
			&\quad+\sum_{i<j}|a_{ij}^{\alpha,\beta}|\delta_h^{i,j}u_h(x).
	\end{split}
\end{equation*}
Under Assumption \ref{assump_dsquare}, one obtains
\begin{equation}\label{eq_fh}
	F_h(D^2_hu)=F(D^2u)+O(h^2),
\end{equation}
provided $u\in C^4(\Omega)$.
}

%

For $0<h<h_0$ fixed, though arbitrary, we denote with $\mathcal{F}_h$ the set of functions defined on $\Omega_h$. A \emph{numerical scheme} (or numerical method, or approximation scheme) is a family $\left(G_h\right)_{h\in(0,h_0)}$ of rules $G_h:\mathcal{F}_h\times\Omega_h\to\mathbb{R}$. We abuse terminology and sometimes refer to a single rule $G_h$ as a numerical scheme. A solution to the numerical scheme is a function $u_h\in\mathcal{F}_h$ such that $G_h(u_h(x),x)=0$ for every $x\in\Omega_h$. Typically, $G_h$ depends on $u_h$ through $D_h^2u_h$, $D_hu_h$ and $u_h$ and we are interested in 
\[
	G_h(u_h(x),x)=G_h(D^2_hu_h(x),D_hu_h(x),u_h(x),x).
\]
Combining the former notation with the discretisations introduced above, we notice $G_h(u_h(x),x)$ also depends on $u_h$ in neighbouring points $y\in\Omega_h$. That is, 
\[
	G_h=G_h\left(u_h(x),\restr{u_h(y)}{y\in N(x)},x\right),
\]
where $N(x)$ stands for the neighbouring points of $x$ used in the required discretisations. We use this notation only in the definition of degenerate elliptic schemes. Elsewhere in the paper, we adhere to $G_h=G_h(u_h(x),x)$ and leave the dependence on the neighbouring points implicit. We continue with the definition of degenerate elliptic methods.
\begin{definition}[Degenerate elliptic]\label{def_degel}
We say the numerical method $(G_h)_{0<h<h_0}$ is degenerate elliptic if
\[
	G_h(v_h(x),x)=G_h\left(v_h(x),\restr{v_h(y)}{y\in N(x)},x\right)
\]
is non-decreasing with respect to $v_h(x)$ and non-increasing with respect to $v_h(y)$, for every $y\in N(x)$.
\end{definition}

Notice that both ${\color{r2}F_h(D^2_hu_h)}$ and $D_hu_h$ are degenerate elliptic discretisations. Next, we recall the definitions of monotonicity, consistency and stability for an approximation scheme.

\begin{definition}[Monotone scheme]\label{def_monotonicity}
We say the numerical method $(G_h)_{0<h<h_0}$ is monotone if, for every $0<h<h_0$, and every $u_h,v_h:\Omega_h\to\mathbb{R}$ with $u_h(x)=v_h(x)$ and $u_h\leq v_h$, we have
\[
	G_h(v_h(x),x)\leq G_h(u_h(x),x).
\]
\end{definition}
The next definition concerns the consistency of the numerical scheme. 
\begin{definition}[Consistent scheme]\label{def_consistency}
We say the numerical method $(G_h)_{0<h<h_0}$ is consistent with \eqref{eq_sundayroast} if
\[
	\limsup_{\substack{h\to 0\\y\to x,\,y\in\Omega_h\\\xi\to 0}}G_h(\varphi(y)+\xi{\color{r1}+e^{-1/h}},y)\leq F^*(D^2\varphi(x), D\varphi(x), \varphi(x),x)
\] 
and
\[
	\liminf_{\substack{h\to 0\\y\to x,\,y\in\Omega_h\\\xi\to 0}}G_h(\varphi(y)+\xi{\color{r1}+e^{-1/h}},y)\geq F_*(D^2\varphi(x), D\varphi(x), \varphi(x),x),
\] 
for every $\varphi\in C^\infty(\overline\Omega)$ and $x\in \overline\Omega$.
\end{definition}

{\color{r1}
\begin{remark}[Attainability of the supremum]\label{rem_teso}\normalfont
In the abstract framework introduced in \cite{BarSou}, driven by the solution operator $S(\rho,x,u^\rho(x),u^\rho)$, the attainability of the supremum and infimum for $u^\rho-\varphi$ may fail. The reason is the fact that $u^\rho$ is merely bounded, but no continuity information is available for the solution to $S=0$. To correct that, the term $\varphi+\xi$ can be replaced by $\varphi+\xi+\eta_\rho$, where $\eta_\rho\to0$, as $\rho\to0$, with a decay related to the approximation structures involved in the analysis. As a consequence, the slack ensures a maximum/minimum is attained by the difference $u^\rho-\varphi$. For example, by choosing $\eta_\rho\sim e^{-1/\rho}$, one ensures the slack decays faster than any power of $\rho$, yielding attainability and preserving the limits involved in the consistency definition.
\end{remark}
}

{\color{wethepeople}
\begin{remark}[Half-relaxed limits through grid points]\label{rem_newconsistency}
In Defintion \ref{def_consistency}, we have considered approximating sequences $(y_n)_{n\in\mathbb{N}}$ in the grid $\Omega_h$. We notice the proof of \cite[Theorem 2.1]{BarSou} can be adjusted to this setting. Indeed, the necessary adjustments involve the use of a realising sequence. Although such a sequence may fail extremality, the slack $e^{-1/h}$ allows us to produce an associated sequence of almost minimisers. Such a sequence converges to the point of interest in the domain and preserves a set of inequalities relevant to the sub and supersolution properties.
\end{remark}
}

We proceed with defining stability. As usual, a function $u_h:\Omega_h\to\mathbb{R}$ denotes a solution to $G_h=0$.
\begin{definition}[Stability]\label{def_stability}
We say the numerical method $(G_h)_{0<h<h_0}$ is stable if there exists $C>0$ such that 
\[
	\sup_{h\in(0,h_0)}\max_{x\in\Omega_h}\left|u_h(x)\right|\leq C.
\]
\end{definition}

Definitions \ref{def_monotonicity}-\ref{def_stability} are the main ingredients in the criterion for convergence of the numerical scheme introduced in \cite{BarSou}. {\color{r2}Indeed, define $\overline u,\underline u:\overline\Omega\to\mathbb{R}$ as
\[
	\overline u(x)\coloneqq \lim_{h\to 0}\sup_{\substack{y\to x\\y\in\Omega_h}}u_h(y)\hspace{.2in}\mbox{and}\hspace{.2in}\underline u(x)\coloneqq \lim_{h\to 0}\inf_{\substack{y\to x\\y\in\Omega_h}}u_h(y)
\]
Under consistency, monotonicity and stability, \cite[Theorem 2.1]{BarSou} ensures
\[
	\begin{cases}
		G(D^2\overline u, D\overline u,\overline u,x)\leq f(x)&\hspace{.2in}\mbox{in}\hspace{.2in}\Omega\\
		G(D^2\underline u, D\underline u,\underline u,x)\geq f(x)&\hspace{.2in}\mbox{in}\hspace{.2in}\Omega.
	\end{cases}
\]
The same result does not enforce the pointwise boundary condition imposed in \eqref{eq_main1} or in \eqref{eq_main2}. To address that issue, it suffices to verify that $\overline u\leq \underline u$ on $\partial\Omega$. Indeed, if that inequality is available and the equation $G=f$ has a comparison principle (in the sense of \cite[Theorem 3.3]{CraIshLio} or, more concretely, in the sense of \cite[Proposition 4]{HuaPimRamSwi}), one concludes $\overline u=\underline u$. Consequently, the convergence of the method follows\footnote{{\color{r2}{The role of the different notions of boundary conditions, and the importance of a comparison principle available for the relaxed Dirichlet problem, was brought to our attention by a anonymous referee. The idea of using a contradiction argument to establish $\overline u\leq g$ and $\underline u\geq g$ on $\partial\Omega$ was also a contribution of the anonymous referee.}}}. 

In the present paper, the convergence of the numerical method is established in two steps. First, consistency, monotonicity and stability ensures that $\overline u$ and $\underline u$ are, respectively, a sub and a supersolution to the limiting equation in $\Omega$. Then we resort to a contradiction argument to establish $\overline u\leq g\geq \underline u$. The comparison principle, as in \cite[Proposition 4]{HuaPimRamSwi}, ensures $\overline u=\underline u$ and complete the proof.
}

%

We close this section with a strategy to solve $G_h=0$. Namely, we consider the solution operator 
\[
	S_\rho[u(x)]\coloneqq u(x)-\rho G_h(u(x),x),
\]
where $0<\rho\ll1$ is a parameter chosen to ensure, among other things, that $S_\rho$ has a fixed point; see \cite{Obe}. The choice of $\rho$ depends on $h$ through a non-linear CFL condition. We also refer to $S_\rho$ as \emph{Euler operator}. We resort to this strategy in the proof of Theorem \ref{thm_convergencepure} as well as in our numerical tests. Our last preliminary result regards the conditions ensuring that the solution operator is a contraction. See \cite[Theorem 7]{Obe}.

\begin{proposition}[The solution operator is a contraction]\label{prop_contraction}
Let $h\in(0,h_0)$ be fixed, though arbitrary. Suppose $G_h$ is degenerate elliptic and Lipschitz-continuous, with Lipschitz constant $C=C(d,\lambda,\Lambda,\varepsilon,h,\theta)$. {\color{r2}Suppose, in addition, that $G_h$ is proper.} If {\color{r2}$0<\rho<C^{-1}$}, the operator $S_\rho$ is a strict contraction.
\end{proposition}

\section{A detour on the existence of solutions}\label{sec_newproofeist}

For completeness, we detail the proof of Theorem \ref{thm_existence}. We argue as usual in the theory of viscosity solutions. Namely, for $\varepsilon>0$ fixed, we prove a comparison principle for \eqref{eq_reg1} and construct global sub and supersolutions. Then Perron's method yields the existence of a viscosity solution $u^\varepsilon$ to that problem. At this point, regularity estimates for \eqref{eq_reg1} allow us to consider the limit $\varepsilon\to 0$. The stability of viscosity solutions yields the existence of a function $u\in C(\Omega)$ solving \eqref{eq_main1}. We proceed with the comparison principle.

\begin{proposition}[Comparison principle]\label{prop_cp}
Suppose Assumptions \ref{assump_A40} and \ref{assump_data} are in force. Let $u^\varepsilon\in{\rm USC}(\overline\Omega)$ be a viscosity subsolution to the PDE in \eqref{eq_reg1} and $w^\varepsilon\in {\rm LSC}(\overline\Omega)$ be a viscosity supersolution to the PDE in \eqref{eq_reg1}. Suppose further that  $u^\varepsilon\leq w^\varepsilon$ on $\partial \Omega$. Then $u^\varepsilon\leq w^\varepsilon$ in $\overline\Omega$.
\end{proposition}
\begin{proof}
The argument follows along the same lines as in the proof of \cite[Proposition 4]{HuaPimRamSwi} and is omitted.
\end{proof}

We proceed by constructing global sub and supersolutions.

\begin{proposition}[Existence of global sub and supersolutions]\label{prop_barriers}
Suppose Assumptions \ref{assump_A40} and \ref{assump_data} are in force. Then, for every $\varepsilon\in(0,1)$, there exists a viscosity subsolution $\underline w\in C(\overline\Omega)$ to \eqref{eq_reg1} and a viscosity supersolution $\overline w\in C(\overline\Omega)$ to \eqref{eq_reg1}. In addition, $\underline w=\overline w=g$ on $\partial\Omega$. Finally, $\underline w$ and $\overline w$ do not depend on $\varepsilon$.
\end{proposition}
\begin{proof}
The argument follows along the same lines as in the proof of \cite[Lemma 2]{HuaPimRamSwi} and is omitted.
\end{proof}

For every $0<\varepsilon<1$,  Perron's method ensures the existence of a viscosity solution $u^\varepsilon\in C(\Omega)$ to \eqref{eq_reg1}, agreeing with $g$ on the boundary of $\Omega$. In addition, $\underline w\leq u^\varepsilon\leq \overline w$, for every $\varepsilon\in (0,1)$. We continue with an observation on the H\"older-regularity of the family $(u^\varepsilon)_{\varepsilon\in (0,1)}$.

\begin{lemma}[Uniform compactness]\label{lem_compactness}
Suppose Assumptions \ref{assump_A40} and \ref{assump_data} are in force.  Then there exists a unique viscosity solution $u^\varepsilon\in C(\Omega)$ to \eqref{eq_reg1}. Moreover, we have $\underline w\leq u\leq \overline w$. Finally, there exists $\alpha\in (0,1)$ such that $u^\varepsilon\in C^\alpha_{\rm loc}(\Omega)$ and, for every $\Omega'\Subset\Omega$, we have
\[
	\left\|u^\varepsilon\right\|_{C^\alpha(\Omega')}\leq C,
\]
where $C=C(d,\lambda,\Lambda,\left\|f\right\|_{L^\infty(\Omega)},\left\|g\right\|_{L^\infty(\partial\Omega)},{\rm dist}(\Omega',\partial\Omega))$.
\end{lemma}
\begin{proof}
Notice that $u^\varepsilon$ satisfies
\[
	-\left\|f\right\|_{L^\infty(\Omega)}-\left\|u^\varepsilon\right\|_{L^\infty(\Omega)}\leq F(D^2u^\varepsilon)\leq \left\|f\right\|_{L^\infty(\Omega)}+\left\|u^\varepsilon\right\|_{L^\infty(\Omega)}
\]
in $\Omega\cap\{|Du|>1\}$. Arguing as in the proof of \cite[Corollary 2]{HuaPimRamSwi}, one obtains the result.
\end{proof}

\begin{proof}[Proof of Theorem \ref{thm_existence}]
The family $(u^\varepsilon)_{\varepsilon\in(0,1)}$ is uniformly bounded in some H\"older space. Therefore, there exists a convergent subsequence, still denoted with $(u^\varepsilon)_{\varepsilon\in(0,1)}$, satisfying $u^\varepsilon\longrightarrow u$, locally uniformly in $\Omega$, where $u\in C(\Omega)$. Standard results on the stability of viscosity solutions ensure that $u$ solves \eqref{eq_main1} in the viscosity sense.
\end{proof}

\section{Degenerate fully nonlinear equations}\label{sec_discrete1}

 We now introduce a discrete approximation of \eqref{eq_reg1}. To ease notation, we drop the superscript $\varepsilon$. The discretisation $G_h$ is given by
\begin{equation}\label{eq_discretea}
	G_h(u_h,x)\coloneqq
		\begin{cases}
			\varepsilon u_h(x)+F_h(D^2_hu_h(x))-\frac{f(x)}{\left(\varepsilon+|D_hu_h(x)|^2\right)^{\frac{\theta}{2}}}&\,\mbox{if}\,x\in\Omega_h\\
			{\color{r2}u_h(x)-g(x)}&\,\mbox{if}\,x\in\partial \Omega_h,
		\end{cases}
\end{equation}

We proceed with the monotonicity of $G_h(u_h,x)$.

\begin{proposition}[Monotonicity of $G_h$]\label{prop_monotone}
Let $G_h$ be defined as in \eqref{eq_discretea}. Suppose Assumptions \ref{assump_dsquare} and \ref{assump_data} are in force. Suppose $u_h,v_h:\Omega_h\to\mathbb{R}$ are such that $u_h(x)=v_h(x)$ for some $x\in \Omega_h$, with $u_h\leq v_h$ in $\Omega_h$. Then $G_h(v_h,x)\leq G_h(u_h,x)$.
\end{proposition}
\begin{proof}
If $x\in\partial\Omega_h$, we conclude 
\[
	{\color{r2}G_h(u_h,x)=u_h(x)-g(x)=v_h(x)-g(x)=G_h(v_h,x)}. 
\]
Suppose otherwise that $x\in\Omega_h$. Under Assumption \ref{assump_dsquare}, we have $F_h(D^2_hv_h(x))\leq F_h(D^2_hu_h(x))$. Also, $|D_hv_h(x)|^2\leq |D_hu_h(x)|^2$. Hence, 
\begin{align*}
		G_h(u_h,x)&=\varepsilon v_h(x)+F_h(D^2_hu_h(x))-\frac{f(x)}{\left(\varepsilon+|D_hu_h(x)|^2\right)^\frac{\theta}{2}}\\
			&\geq \varepsilon v_h(x)+F_h(D^2_hv_h(x))-\frac{f(x)}{\left(\varepsilon+|D_hv_h(x)|^2\right)^\frac{\theta}{2}}\\
			&=G_h(v_h,x),
\end{align*}
which completes the proof.
\end{proof}

\begin{remark}[The case of $p$-Laplace type equations]\label{rem_espen1} 
The $p$-Laplace and the porous medium equations are natural counterparts to \eqref{eq_main1} in the divergence-form setting. In their cases, the product of monotone operators may also lead to a scheme lacking monotonicity. We believe our strategy can be adjusted to address those cases as well. For a monotone discretisation of the porous medium equation in dimension $d=1$, we refer to \cite{delJak}.
\end{remark}

In what follows, we verify that $G_h(u_h,x)$ is Lipschitz-continuous.

\begin{proposition}[Lipschitz continuity of $G_h$]\label{prop_lipcon}
Let $G_h$ be defined as in \eqref{eq_discretea}. Suppose Assumptions \ref{assump_dsquare} and \ref{assump_data} hold true. Then
\[
	\max_{x\in\Omega_h}\left|G_h(u_h,x)-G_h(v_h,x)\right|\leq C\max_{x\in\Omega_h}\left|u_h(x)-v_h(x)\right|,
\]
for every $u_h,v_h:\Omega_h\to\mathbb{R}$. Moreover, the constant $C>0$ depends on the dimension $d$, the ellipticity constants $0<\lambda\leq\Lambda$, $0<\varepsilon\ll1$, $0<h\ll1$, and $\left\|f\right\|_{L^\infty(\Omega_h)}$. 
\end{proposition}
\begin{proof}
If $x\in \partial\Omega_h$, $G_h(u_h,x)-G_h(v_h,x)={\color{r2}u_h(x)-v_h(x)}$ and there is nothing to prove. Suppose otherwise and notice
\begin{equation}\label{eq_wow}
	\begin{split}
		\left|G_h(u_h,x)-G_h(v_h,x)\right|&\leq \left(\varepsilon+\frac{C(\lambda,\Lambda,d)}{h^2}\right)\left|u_h(x)-v_h(x)\right|\\
			&\quad+\left|\frac{f(x)}{\left(\varepsilon+|D_hu_h(x)|^2\right)^\frac{\theta}{2}}-\frac{f(x)}{\left(\varepsilon+|D_hv_h(x)|^2\right)^\frac{\theta}{2}}\right|\\
			&\eqqcolon \left(\varepsilon+\frac{C(\lambda,\Lambda,d)}{h^2}\right)\left|u_h(x)-v_h(x)\right|+I.
	\end{split}
\end{equation}
We proceed by examining the term $I$. Indeed,
\begin{equation}\label{eq_abomination}
	\begin{split}
		I&\leq\left\|f\right\|_{L^\infty(\Omega)}\left|\frac{\left(\varepsilon+|D_hu_h(x)|^2\right)^\frac{\theta}{2}-\left(\varepsilon+|D_hv_h(x)|^2\right)^\frac{\theta}{2}}{\left(\varepsilon+|D_hu_h(x)|^2\right)^\frac{\theta}{2}\left(\varepsilon+|D_hv_h(x)^2\right)^\frac{\theta}{2}}\right|
	\end{split}
\end{equation}
and
\begin{equation}\label{eq_wow2}
	\begin{split}
		&\left|\frac{\left(\varepsilon+|D_hu_h(x)|^2\right)^\frac{\theta}{2}-\left(\varepsilon+|D_hv_h(x)|^2\right)^\frac{\theta}{2}}{\left(\varepsilon+|D_hu_h(x)|^2\right)^\frac{\theta}{2}\left(\varepsilon+|D_hv_h(x)^2\right)^\frac{\theta}{2}}\right|\\
			&\leq\frac{\max\left(\left(\varepsilon+|D_hu_h(x)|^2\right),\left(\varepsilon+|D_hv_h(x)|^2\right)\right)^{\frac{\theta}{2}-1}}{\left(\varepsilon+|D_hu_h(x)|^2\right)^\frac{\theta}{2}\left(\varepsilon+|D_hv_h(x)|^2\right)^\frac{\theta}{2}}\left||D_hu_h(x)|^2-|D_hv_h(x)|^2\right|\\
			&\leq \frac{C(d,\varepsilon)}{h}\max_{x\in\Omega_h}\left|u_h(x)-v_h(x)\right|.
	\end{split}
\end{equation}

By combining \eqref{eq_wow} and \eqref{eq_abomination}, and \eqref{eq_wow2}, we conclude
\[
	\max_{x\in\Omega_h}\left|G_h(u_h(x),x)-G_h(v_h(x),x)\right|\leq C\max_{x\in\Omega_h}\left|u(x)-v_h(x)\right|,
\]
where
\begin{equation}\label{eq_mrwa}
	C\coloneqq\left(\varepsilon+\frac{C(\lambda,\Lambda,d)}{h^2}+\frac{C(d,\varepsilon)}{h}\right)\sim\frac{1}{h^2}.
\end{equation}
\end{proof}

\begin{remark}[Courant-Friedrichs-Lewy condition]\label{rem_cfl}
In the sequel, we resort to an Euler scheme with parameter $0<\rho\ll1$ to prove the existence of solutions to $G_h=0$. In that context, the constant in \eqref{eq_mrwa} is pivotal. Indeed, we impose a Courant-Friedrichs-Lewy condition of the form
\[
	\rho< \left(\varepsilon+\frac{C(\lambda,\Lambda,d)}{h^2}+\frac{C(d,\varepsilon)}{h}\right)^{-1}.
\]
Having in mind that our goal is to examine the case $0<\varepsilon,h\ll1$, we can simplify the former inequality and work under {\color{r2}$\rho\ll h^2$}.
\end{remark}

We continue with the consistency of $G_h$. To that end, write 
\begin{equation}\label{eq_reg12}
	G(D^2u,Du,u,x)\coloneqq
		\begin{cases}
			\varepsilon u(x)+F(D^2u(x))-\frac{f(x)}{\left(\varepsilon+|Du(x)|^2\right)^\frac{\theta}{2}}&\hspace{.1in}\mbox{in}\hspace{.1in}\Omega\\
			u(x)-g(x)&\hspace{.1in}\mbox{on}\hspace{.1in}\partial\Omega.
	\end{cases}
\end{equation}
A viscosity solution to $G=0$ solves \eqref{eq_reg1} in the viscosity sense. We prove that $G_h$ is consistent with $G$.

\begin{proposition}[Consistency]\label{prop_consistency}
Suppose {\color{r2}Assumptions \ref{assump_dsquare}-\ref{assump_data} hold}. Then $G_h$ is consistent with $G$.
\end{proposition}
\begin{proof}
For $\varphi\in C^\infty(\overline\Omega)$ and $x\in\overline\Omega$, we prove that 
\begin{equation}\label{eq_cons1}
	\limsup_{\substack{h\to 0\\y\to x,\,y\in\Omega_h\\\xi\to 0}}G_h(\varphi(y)+\xi+e^{-1/h},y)\leq G^*\left(D^2\varphi(x),D\varphi(x),\varphi(x),x\right).
\end{equation}
We split the proof into two steps, depending on the point $x\in\overline\Omega$. We start by considering $x\in\Omega$.

\noindent{\bf Step 1 - }If $x\in\Omega$, we can suppose the points $y\in\Omega_h$ approaching it are also interior points. In this case,
\begin{align*}
		G_h(\varphi(y)+\xi+e^{-1/h},y)&=\varepsilon (\varphi(y)+\xi+e^{-1/h})+F_h(D^2_h(\varphi(y)+\xi+e^{-1/h}))\\&\quad-\frac{f(y)}{\left(\varepsilon+|D_h(\varphi(y)+\xi+e^{-1/h})|^2\right)^\frac{\theta}{2}}
\end{align*}

The regularity of $\varphi$ implies
\[
	D_h(\varphi(y)+\xi+e^{-1/h})= D\varphi(x)+O(h).
\]
In addition, \eqref{eq_fh} builds upon the regularity of $\varphi$ to ensure
\begin{align*}
		\limsup_{\substack{h\to 0\\y\to x,\,y\in\Omega_h\\\xi\to 0}}G_h(\varphi(y)+\xi+e^{-1/h},y)&=\varepsilon \varphi(x)+F(D^2\varphi(x))-\frac{f(x)}{\left(\varepsilon+|D\varphi(x)|^2\right)^\frac{\theta}{2}}\\
			&\leq G^*\left(D^2\varphi(x),D\varphi(x),\varphi(x),x\right).
\end{align*}

\noindent{\bf Step 2 - }Consider next $x\in\partial\Omega$. In this case, one can approach $x$ by points $y\in\Omega_h$ or $y\in\partial\Omega_h$, or both. Since we work with limit superiors, we consider only $y\in\Omega_h$ or $y\in\partial\Omega_h$. In the latter case, we have
\[
	G_h(\varphi(y)+\xi+e^{-1/h},y)=\varphi(y)+\xi+e^{-1/h}-g(y)\longrightarrow\varphi(x)-g(x)=G(D^2\varphi,D\varphi,\varphi,x)
\]
as $h\to0$, $y\to x$ and $\xi\to 0$, and \eqref{eq_cons1} follows. Suppose now $y\in\Omega_h$. Then, arguing as before,
\begin{align*}
		&\limsup_{\substack{h\to 0\\y\to x,\,y\in\Omega_h\\\xi\to 0}}G_h(\varphi_{\xi,h}(y),y)\\
			&\leq \limsup_{\substack{h\to 0\\y\to x\\\xi\to 0}}\left(\varepsilon \varphi_{\xi,h}(y)+F_h(D^2_h\varphi_{\xi,h}(y))-\frac{f(y)}{\left(\varepsilon+|D_h\varphi_{\xi,h}(y)|^2\right)^\frac{\theta}{2}}\right)\\
			&\leq G^*\left(D^2\varphi(x),D\varphi(x),\varphi(x),x\right),
\end{align*}
where
\[
	\varphi_{\xi,h}(y)\coloneqq\varphi(y)+\xi+e^{-1/h},
\]
\end{proof}

The consistency of the method is fundamental for convergence. However, it also plays a role in building global barriers for $G_h$. Indeed, by modifying the sub and supersolutions in Proposition \ref{prop_barriers}, one can find $\underline w$ and $\overline w$ such that $\underline w=\overline w=g$ on $\partial\Omega_h$ and $G_h(\underline w(x),x)\leq 0\leq G_h(\overline w(x),x)$ in $\Omega_h$ for every $0<h\ll1$ small enough. We formalise this heuristic in the next proposition.

\begin{proposition}[Discrete global barriers]\label{prop_discretegb}
Suppose {\color{r2}Assumptions \ref{assump_dsquare}-\ref{assump_data} hold true}. There exists $0<h_0\ll1$ such that, for $0<h<h_0$, one can find $\underline w,\overline w:\Omega_h\to\mathbb{R}$ with $\underline w\leq g\leq \overline w$ on $\partial\Omega_h$, satisfying
\[
	G_h(\underline w(x),x)\leq 0\leq G_h(\overline w(x),x),
\]
for every $x\in \Omega_h$.
\end{proposition}
\begin{proof}
For constants $C_1,C_2>0$, define 
\[
	\overline w(x)\coloneqq C_2-\frac{C_1}{2\lambda d}\left|x-x_0\right|^2,
\]
where $x_0\in\mathbb{R}^d$ is such that $|x_i-x_{0,i}|>\sqrt{\lambda}$ for every $i=1,\ldots,d$, and $C_2$ is such that $\overline w(x)\geq \left\|g\right\|_{L^\infty(\partial\Omega)}$. Note
\[
	D^2_h\overline w(x)={\color{r2}-\frac{C_1}{\lambda d}I},
\]
whereas
\[
	|D_h{\color{r2}\overline{w}}(x)|^2=\frac{C_1^2}{4(\lambda d)^2}\sum_{i=1}^d\left(\max\left(2|x_i-x_{0,i}|-h,0\right)\right)^2.
\]
Now, let $x\in\partial\Omega_h$. Then $G_h(\overline w(x),x)=\overline w(x)-g(x)\geq 0$, because of the choice of $C_1$. If $x\in \Omega_h$, we have
\[
	G_h(\overline w(x),x)\geq \varepsilon\overline w(x)+C_1-\frac{C_1}{(\varepsilon C_1^2)^\frac{\theta}{2}}\geq 0.
\]
Hence, $G_h(\overline w(x),x)\geq 0$ for every $x\in\Omega_h$. The construction of $\underline w$ is entirely analogous.
\end{proof}

 An important aspect of Proposition \ref{prop_discretegb} concerns uniform bounds on $\underline w$ and $\overline w$. Because these functions are obtained as (uniform) variants of the global barriers in Proposition \ref{prop_barriers}, we conclude there exists $C>0$, depending only on the dimension $d$, ellipticity $0<\lambda\leq \Lambda$, and the norms $\left\|f\right\|_{L^\infty(\Omega)}$ and $\left\|g\right\|_{L^\infty(\partial\Omega)}$ such that
 \[
 	-C\leq \underline w\leq \overline w\leq C.
 \]
 It is critical to notice that $C$ does not depend on $h$.
 
\begin{proposition}[Stability]\label{prop_stability}
Let $h\in(0,h_0)$ be fixed. Let $u_h:\Omega_h\to\mathbb{R}$ be a solution to $G_h(u_h(x),x)=0$ in $\Omega_h$. Suppose {\color{r2}Assumptions \ref{assump_dsquare}-\ref{assump_data} are in force}. Then there exists $C>0$ such that 
\[
	\sup_{h\in(0,h_0)}\max_{x\in\Omega_h}\left|u_h(x)\right|\leq C.
\]
In addition, $C$ depends on the dimension $d$, the ellipticity $0<\lambda\leq \Lambda$, and the norms $\left\|f\right\|_{L^\infty(\Omega)}$ and $\left\|g\right\|_{L^\infty(\partial\Omega)}$, but does not depend on $h$.
\end{proposition}
\begin{proof}
Let $\overline w$, be the barrier function from Proposition \ref{prop_discretegb}. We claim that $u_h\leq \overline w$ in $\overline \Omega_h$. Suppose otherwise; if this is the case, there exists $\overline x\in \Omega_h$ such that $u_h(\overline x)>\overline w(\overline x)$. Clearly, such a point has to be in $\Omega_h$. Also, $u_h(\overline x)\geq u_h(x)$ for every $x\in \Omega_h$. Hence,
\[
	u_h(\overline x)-u_h(x)\geq \overline w(\overline x)-\overline w(x),
\]
for every $x\in \Omega_h$. It follows that 
\begin{equation*}\label{eq_itsthend}
	\begin{split}
		G_h(u_h(\overline x),\overline x)&=\varepsilon u_h(\overline x){\color{r2}+}F_h(D^2_hu_h(x))-\frac{f(\overline x)}{\left(\varepsilon+|D_hu_h(\overline x)|^2\right)^\frac{\theta}{2}}\\
			&>\varepsilon \overline w(\overline x){\color{r2}+}F_h(D^2_h\overline w(x))-\frac{f(\overline x)}{\left(\varepsilon+|D_h\overline w(\overline x)|^2\right)^\frac{\theta}{2}}.
	\end{split}
\end{equation*}
Therefore, $G_h(\overline w(\overline x),\overline x)<G_h(u_{\color{r2}h}(\overline x),\overline x)=0$, which contradicts Proposition \ref{prop_discretegb}. We conclude that $u_h\leq \overline w$. Arguing as before, with $\underline w$ instead of $\overline w$, one obtains
\[
	u_h\geq -\overline w
\]
in $\overline \Omega_h$ and completes the proof.
\end{proof}

We have established that the method in \eqref{eq_discretea} is monotone, consistent with \eqref{eq_reg1} and stable. Therefore, we are in a position to prove Theorem \ref{thm_convergencepure}.

\begin{proof}[Proof of Theorem \ref{thm_convergencepure}]
{\color{r2}Define $\overline u,\underline u:\overline\Omega\to\mathbb{R}$ as
\[
	\overline u(x)\coloneqq\lim_{h\to 0}\sup_{\substack{y\to x\\y\in\Omega_h}}u_h^\varepsilon(y)\hspace{.2in}\mbox{and}\hspace{.2in}\underline u(x)\coloneqq\lim_{h\to 0}\inf_{\substack{y\to x\\y\in\Omega_h}}u_h^\varepsilon(y).
\]
For ease of presentation, we split the proof into three steps.

\medskip

\noindent{\bf Step 1 - }Because the method $G_h$ is monotone, stable and consistent with \eqref{eq_reg1}, we infer from \cite[Theorem 2.1]{BarSou}, that
\begin{equation}\label{eq_solintp}
		\begin{cases}
			\varepsilon \overline u(x)+F(D^2\overline u(x))\leq \dfrac{f(x)}{\left(\varepsilon+|D\overline u(x)|^2\right)^\frac{\theta}{2}}&\hspace{.1in}\mbox{in}\hspace{.1in}\Omega\\
			\varepsilon \underline u(x)+F(D^2\underline u(x))\geq \dfrac{f(x)}{\left(\varepsilon+|D\underline u(x)|^2\right)^\frac{\theta}{2}}&\hspace{.1in}\mbox{in}\hspace{.1in}\Omega.
		\end{cases}
\end{equation}
To conclude the convergence of the method $G_h$ as $h\to 0$, it remains to examine the boundary behaviour of $\overline u$ and $\underline u$. Suppose that $\overline u\leq \underline u$ on $\partial\Omega$. In that case, \eqref{eq_solintp} builds upon the comparison principle in Proposition \ref{prop_cp} to yield $\overline u\leq \underline u$ in $\Omega$. Naturally, one would have $\overline u=\underline u$ and the convergence of the method. The next step verifies that $\overline u\leq \underline u$ on $\partial\Omega$.

\medskip

\noindent{\bf Step 2 - }The claim is tantamount to $\overline u \leq g\leq \underline u$ on the boundary of $\Omega$. We detail the case $g\geq \overline u$, as the remaining one is entirely analogous. Suppose, by contradiction, that the inequality is false. Then there exists $y\in\partial\Omega$ such that $g(y)>\overline u(y)$. Because $\Omega$ satisfies a uniform exterior condition, there exist $R>0$ and $x_y\in\mathbb{R}^d\setminus\Omega$ such that $R=|y-x_y|$ and $\overline B_R(x_y)\cap\overline\Omega=\{y\}$.

Set 
\[
	\alpha>\max\left(\frac{\Lambda-2\lambda}{\lambda},2\right),
\]
where $0<\lambda\leq\Lambda$ are the ellipticity constants in Assumption \ref{assump_A40}. Define the barrier $\phi_y\in C^2(\overline\Omega)$ as
\[
	\phi_y(x)\coloneqq R^{-\alpha}-|x-x_y|^{-\alpha}.
\]
It is clear that $\phi_y\geq 0$ in $\overline\Omega$, with $\phi_y(x)=0$ if and only if $x=y$. For $M>0$, yet to be chosen, denote by $x_M$ the point
\[
	x_M\coloneqq\arg\max_{x\in\overline\Omega}\overline u(x)-M\phi_y(x).
\]
Since $\overline u$ is an upper semi-continuous function, the point $x_M$ is well defined. The subsolution property for $\overline u$ implies
\begin{equation}\label{eq_subsolp1}
	\varepsilon M\phi_y(x_M)+F(D^2M\phi_y(x_M))\leq \dfrac{f(x_M)}{\left(\varepsilon+|DM\phi_y(x_M)|^2\right)^\frac{\theta}{2}},\hspace{.2in}\mbox{if}\hspace{.1in}x_M\in\Omega,
\end{equation}
and
\begin{equation}\label{eq_subsolp2}
	\min\left(\varepsilon M\phi_y(x_M)+F(D^2M\phi_y(x_M))- \dfrac{f(x_M)}{\left(\varepsilon+|DM\phi_y(x_M)|^2\right)^\frac{\theta}{2}}, \overline u(x_M)-g(x_M)\right)\leq 0,
\end{equation}
if $x_M\in\partial\Omega$. 

Now, the definition of $x_M$ implies
\[
	\lim_{M\to\infty}x_M=y\hspace{.2in}\mbox{and}\hspace{.2in}\lim_{M\to\infty}M\phi_y(x_M)=0.
\]
As a consequence, for $M\gg1$ large enough, we get $\overline u(x_M)-g(x_M)>0$ in case $x_M\in\partial\Omega$. This fact, combined with \eqref{eq_subsolp1} and \eqref{eq_subsolp2}, yields
\begin{equation}\label{eq_subsolp3}
	\varepsilon M\phi_y(x_M)+F(D^2M\phi_y(x_M))\leq \dfrac{f(x_M)}{\left(\varepsilon+|DM\phi_y(x_M)|^2\right)^\frac{\theta}{2}},\hspace{.2in}\mbox{for}\hspace{.1in}x_M\in\overline\Omega.
\end{equation}
However,
\[
	\varepsilon M\phi_y(x_M)+F(D^2M\phi_y(x_M))\geq \frac{\alpha M}{|x-x_M|^{\alpha+2}}\left(\lambda(\alpha+2)-\Lambda\right).
\]
The definition of $\alpha>2$ allows us to take $M\gg1$ large enough so that
\[
	\varepsilon M\phi_y(x_M)+F(D^2M\phi_y(x_M)) >\dfrac{\left\|f\right\|_{L^\infty(\Omega)}}{\varepsilon^\frac{\theta}{2}},
\]
obtaining a contradiction with \eqref{eq_subsolp3} and concluding that $\overline u\leq g$ on $\partial\Omega$.

\medskip

\noindent{\bf Step 3 - }Finally, because $G_h$ is degenerate elliptic and Lipschitz continuous, with a Lipschitz constant of the order $h^{-2}$, the solution operator $S_\rho$ is a strict contraction, provided we choose $\rho\ll h^2$. Therefore, seen as a functional on $\ell_\infty$, $S_\rho$ admits a unique fixed point. To complete the proof, we recall $u^\varepsilon$ converges locally uniformly, through a subsequence if necessary, to a viscosity solution to \eqref{eq_main1}.}
\end{proof}

\section{The fully nonlinear free boundary problem}\label{sec_pure}

{\color{r2}We start by describing the construction of an auxiliary function. Let $v\in C(\overline\Omega)$ and define $g^v_\varepsilon:\overline\Omega\to\mathbb{R}$ as 
\[
	g_\varepsilon^v(x)\coloneqq\max\left(\min\left(\frac{v+\varepsilon}{2\varepsilon},1\right),0\right).
\]
Let $h_\varepsilon^v(x)\coloneqq\left(g_\varepsilon^v\ast\eta_\varepsilon\right)(x)$, where $\eta_\varepsilon$ is a standard mollifier, and define
\[
	\theta_\varepsilon^v(x)\coloneqq h_\varepsilon^v(x)\theta_1+(1-h_\varepsilon^v(x))\theta_2.
\]
}
{\color{r2}

We examine the following doubly-regularised problem
\begin{equation}\label{eq_reg2b}
	\begin{cases}
		\varepsilon u^{\varepsilon,v}+F(D^2u^{\varepsilon,v})-\dfrac{f}{\left(\varepsilon+|Du^{\varepsilon,v}|^2\right)^\frac{\theta^v_\varepsilon(x)}{2}}=0&\hspace{.2in}\mbox{in}\hspace{.1in}\Omega\\
		u^{\varepsilon,v}=g&\hspace{.2in}\mbox{on}\hspace{.1in}\partial\Omega,
	\end{cases}
\end{equation}
and propose the method $G_h^{\varepsilon,v}=G_h^{\varepsilon,v}(u^{\varepsilon,v}_h(x),x)$ defined as
\begin{equation}\label{eq_method2b}
	G_h^{\varepsilon,v}\coloneqq
		\begin{cases}
			\varepsilon u_h^{\varepsilon,v}(x)+F_h(D^2_hu^{\varepsilon,v}_h(x))-\dfrac{f(x)}{\left(\varepsilon + |D_hu_h^{\varepsilon,v}(x)|^2\right)^\frac{\theta_\varepsilon^v(x)}{2}}=0&\hspace{.2in}\mbox{in}\hspace{.1in}\Omega_h\\
			u^{\varepsilon,v}_h(x)=g(x)&\hspace{.2in}\mbox{on}\hspace{.1in}\partial\Omega_h.
		\end{cases}
\end{equation}

To ease notation, we drop the superscript $\varepsilon$ in what follows and write $u_h^{\varepsilon,v}=u_h^v$ and $G_h^{\varepsilon,v}=G_h^v$.

{\color{r2}
\begin{remark}[Comparison principle, existence and uniqueness of solutions and regularity estimates]\label{rem_pr}
Our analysis depends on the availability of a comparison principle for \eqref{eq_reg2b}, as well as on the existence and uniqueness of solutions, together with regularity estimates yielding compactness. We notice the comparison principle was established in \cite[Proposition 4]{HuaPimRamSwi}, whereas the existence and uniqueness of solutions was proved in \cite[Theorem 1]{HuaPimRamSwi}. The compactness stems from regularity estimates in H\"older spaces obtained in Corollary \cite[Corollary 1]{HuaPimRamSwi}.
\end{remark}
}

Now, we verify that $G_h^v$ is monotone in the sense of Definition \ref{def_monotonicity}.
\begin{proposition}[Monotonicity of $G_h^v$]\label{prop_monotone2}
Let $G_h^v$ be defined as in \eqref{eq_method2b}. Suppose Assumptions \ref{assump_dsquare} and \ref{assump_data} are in force. Suppose $u_h,w_h:\Omega_h\to\mathbb{R}$ are grid functions such that $u_h(x)=w_h(x)$ for some $x\in \Omega_h$, with $u_h\leq w_h$ in $\Omega_h$. Then $G_h^v(w_h,x)\leq G^v_h(u_h,x)$.
\end{proposition}
\begin{proof}
The proof follows along the same lines as in Proposition \ref{prop_monotone}, noticing that $\theta_\varepsilon^v(x)$ does not depend either on $u_h$ or on $w_h$.
\end{proof}

As before, we verify that $G^v_h$ is consistent. 

\begin{proposition}[Consistency]\label{prop_consistency2}
Suppose Assumptions \ref{assump_dsquare}-\ref{assump_theta} are in force. Then $G_h^v$ is consistent with \eqref{eq_reg2b}.
\end{proposition}
\begin{proof}
The argument is the same as in the proof of Proposition \ref{prop_consistency}, once we notice $\theta_\varepsilon^v\in C(\overline\Omega)$.
\end{proof}

We proceed with the stability of the method. To that end, we build sub and supersolutions to $G^v_h=0$. Once those functions are available, we compare them with solutions $u_h$ at maximum points and take advantage of the asymptotic behaviour of $\theta_\varepsilon$.
\begin{proposition}[Stability]\label{prop_stability2}
Let $h\in(0,h_0)$ be fixed. Let $u^v_h:\Omega_h\to\mathbb{R}$ be a solution to $G^v_h(u^v_h(x),x)=0$ in $\Omega_h$. Suppose assumptions \ref{assump_dsquare}-\ref{assump_theta} are in force. Then
\[
	\left|u^v_h(x)\right|\leq \frac{1}{\varepsilon}\left(\left\|g\right\|_{L^\infty(\partial\Omega)}+1+\frac{\left\|f\right\|_{L^\infty(\Omega)}}{\varepsilon^\frac{\theta_2}{2}}\right).
\]
\end{proposition}
\begin{proof}{\color{r2}
For ease of presentation, we split the proof into three steps. Set $\overline w:\Omega_h\to\mathbb{R}$ as 
\[
	\overline w(x)\coloneqq \frac{1}{\varepsilon}\left(\left\|g\right\|_{L^\infty(\partial\Omega)}+1+\frac{\left\|f\right\|_{L^\infty(\Omega)}}{\varepsilon^\frac{\theta_2}{2}}\right).
\]

\noindent{\bf Step 1 - }We start by verifying that 
\begin{equation}\label{eq_ghw0}
	G^v_h(\overline w(x),x)\geq 0,
\end{equation}	
for every $x\in\Omega_h$. If $x\in\partial\Omega_h$,
\[
	G^v_h(\overline w(x),x)=\overline w(x)-g(x)\geq \frac{1}{\varepsilon}\left(\left\|g\right\|_{L^\infty(\partial\Omega)}+1+\frac{\left\|f\right\|_{L^\infty(\Omega)}}{\varepsilon^\frac{\theta_2}{2}}\right)-\left\|g\right\|_{L^\infty(\partial\Omega_h)}\geq 0.
\]
Now, suppose $x\in\Omega_h$. In this case,
\[
	G^v_h(\overline w(x),x)=\left\|g\right\|_{L^\infty(\partial\Omega)}+1+\frac{\left\|f\right\|_{L^\infty(\Omega)}}{\varepsilon^\frac{\theta_2}{2}}-\frac{\left\|f\right\|_{L^\infty(\Omega)}}{\varepsilon^\frac{\theta_\varepsilon^v(x)}{2}}\geq 0.
\]

\noindent{\bf Step 2 - } We claim that $u^v_h\leq \overline w$ in $\overline \Omega_h$. Suppose otherwise; if this is the case, there exists $\overline x\in \Omega_h$ such that $u^v_h(\overline x)>\overline w(\overline x)$. Such a point has to be in $\Omega_h$. Also, $u^v_h(\overline x)\geq u_h^v(x)$ for every $x\in \Omega_h$. Hence,
\[
	u^v_h(\overline x)-u^v_h(x)\geq \overline w(\overline x)-\overline w(x),
\]
for every $x\in \Omega_h$. It follows that 
\begin{equation}\label{eq_itsthenda}
	\begin{split}
		G^v_h(u^v_h(\overline x),\overline x)&=\varepsilon u^v_h(\overline x){\color{r2}+}F_h(D^2_hu^v_h(x))-\frac{f(\overline x)}{\left(\varepsilon+|D_hu^v_h(\overline x)|^2\right)^\frac{\theta_\varepsilon^v(\overline x)}{2}}\\
			&>\varepsilon \overline w(\overline x){\color{r2}+}F_h(D^2_h\overline w(x))-\frac{f(\overline x)}{\left(\varepsilon+|D_h\overline w(\overline x)|^2\right)^\frac{\theta_\varepsilon^v(\overline x)}{2}}\\
		 	&=G^v_h(\overline  w(\overline x),\overline x).
	\end{split}
\end{equation}
The former observation builds upon \eqref{eq_itsthenda} and the fact that $u^v_h$ is a solution to ensure
\[
	G^v_h(\overline w(\overline x),\overline x)<G^v_h(u_h^v(\overline x),\overline x)=0,
\]
which contradicts \eqref{eq_ghw0}. Therefore, $u^v_h\leq \overline w$.

\noindent{\bf Step 3 - }Consider now $\underline w$ defined as
\[
	\underline w\coloneqq -\frac{1}{\varepsilon}\left(\left\|g\right\|_{L^\infty(\partial\Omega)}+1+\frac{\left\|f\right\|_{L^\infty(\Omega)}}{\varepsilon^\frac{\theta_1}{2}}\right).
\]
Arguing as before, one concludes that
\[
	u^v_h\geq -\overline w
\]
in $\overline \Omega_h$. Hence
\[
	- \frac{1}{\varepsilon}\left(\left\|g\right\|_{L^\infty(\partial\Omega)}+1+\frac{\left\|f\right\|_{L^\infty(\Omega)}}{\varepsilon^\frac{\theta_1}{2}}\right)\leq u^v_h\leq  \frac{1}{\varepsilon}\left(\left\|g\right\|_{L^\infty(\partial\Omega)}+1+\frac{\left\|f\right\|_{L^\infty(\Omega)}}{\varepsilon^\frac{\theta_2}{2}}\right).
\]
Since $\theta_1<\theta_2$ and $0<\varepsilon<1$, it is clear that 
\[
	-\frac{\left\|f\right\|_{L^\infty(\Omega)}}{\varepsilon^\frac{\theta_2}{2}}\leq -\frac{\left\|f\right\|_{L^\infty(\Omega)}}{\varepsilon^\frac{\theta_1}{2}}.
\]
Because the bounds above are independent of $h$, the result follows.}
\end{proof}

We have established that the method in \eqref{eq_method2b} is monotone, consistent with \eqref{eq_reg2b} and stable \emph{independently of $v$ and $\varepsilon$}. Now, we perform a fixed-point analysis. Define the set $\mathcal{B}_h\subset\mathbb{R}^{N(h)}$ as
\begin{equation}\label{eq_bh}
	\mathcal{B}_h\coloneqq\left\lbrace w\in \mathbb{R}^{N(h)}\,|\,\min_{i=1,\ldots,N(h)}\underline w(x_i)\leq w\leq \max_{i=1,\ldots,N(h)}\overline w(x_i)\right\rbrace	.
\end{equation}
It shall be clear that $\mathcal{B}_h$ is closed and convex. Define the map $T_h:\mathcal{B}_h\to\mathbb{R}^{N(h)}$ as follows. For $v\in \mathcal{B}_h$, let $T_hv$ be the unique solution to $G_h^{v}=0$ in $\Omega_h$. We proceed with the study of the map $T_h$.

\begin{proposition}[Properties of the map $T_h$]\label{prop_T}
Fix $0<h\ll1$ and define $\mathcal{B}_h$ as in \eqref{eq_bh}. Then $T_h(\mathcal{B}_h)\subset \mathcal{B}_h$. Also, the map $T_h$ is continuous and precompact.
\end{proposition}
\begin{proof}
We split the argument into three steps. First, we verify $T_h(\mathcal{B}_h)\subset\mathcal{B}_h$.

\smallskip

\noindent{\bf Step 1 - }Let $v\in\mathcal{B}_h$. Then $Tv$ is the unique solution to $G_{h}^v=0$. Because of Proposition \ref{prop_stability2}, we have $\underline w(x)\leq T_hv(x)\leq \overline w(x)$, for every $x\in \Omega_h$. Hence, 
\[
	\min_{j=1,\ldots,N(h)}\underline w(x_j)\leq T_hv(x_i)\leq  \max_{j=1,\ldots,N(h)}\overline w(x_j)
\]
for every $i=1,\ldots,N(h)$. Therefore, $T_hv\in \mathcal{B}_h$.

\smallskip

\noindent{\bf Step 2 - }To verify that $T_h$ is continuous, fix $(v_n)_{n\in\mathbb{N}}\subset \mathcal{B}_h$. Suppose there exists $v_\infty\in \mathcal{B}_h$ such that $v_n\to v_\infty$, as $n\to\infty$. We must verify that $T_hv_n\to T_hv_\infty$, as $n\to\infty$.

In fact, denote with $u_{n}$ the unique solution to $G_{h}^{v_n}=0$. The sequence $(u_n)_{n\in\mathbb{N}}$ is uniformly bounded, because of Proposition \ref{prop_stability2}. Compactness in $\mathbb{R}^{N(h)}$ ensures the existence of a subsequence $(u_{n_k})_{k\in\mathbb{N}}$, converging to a function $\overline u\in\mathcal{B}_h$. Let $x\in \overline \Omega_h$ be fixed, though arbitrary. If $x\in \partial\Omega_h$, we have
\begin{equation}\label{eq_conT1}
	0=G_{h}^{v_{n_k}}(u_{n_k}(x),x)=u_{n_k}(x)-g(x)\to\overline u(x)-g(x).
\end{equation}
If $x\in\Omega_h$, we have 
\begin{equation}\label{eq_conT2}
	0=G_{h}^{v_{n_k}}(u_{n_k}(x),x)\to \varepsilon \overline u(x)+F_h(D^2_h\overline u(x))-\frac{f(x)}{\left(\varepsilon+|D_h\overline u(x)|^2\right)^\frac{\theta_\varepsilon^{v_\infty}(x)}{2}}
\end{equation}
Combining \eqref{eq_conT1} and \eqref{eq_conT2}, one obtains
\[
	G_{h}^{v_\infty}(\overline u(x),x)=0
\]
for every $x\in \Omega_h$. The uniqueness of solutions to $G_h^{v_\infty}=0$ implies that 
\[
	T_hv_\infty=\overline u=\lim_{n\to\infty}T_hv_n,
\]
independently of the subsequence $(u_{n_k})_{k\in\mathbb{N}}$, ensuring the continuity of $T_h$. 

\smallskip

\noindent{\bf Step 3 - }In the sequel, we prove that $T_h$ is precompact. Let $(T_hv_n)_{n\in\mathbb{N}}\subset T(\mathcal{B}_h)$. As before, the stability of the method guarantees that $(Tv_n)_{n\in\mathbb{N}}$ is a uniformly bounded subset of $\mathbb{R}^{N(h)}$. Hence, it admits a convergent subsequence, and we finish the proof.
\end{proof}

The following corollary to the Proposition \ref{prop_T} accounts for the information we use next on the map $T_h$.

\begin{corollary}\label{cor_pf}
Let $0<h\ll1$ be fixed and define $\mathcal{B}_h$ as in \eqref{eq_bh}. There exists $u_h\in \mathcal{B}_h$ such that 
\[
	G_{h}^{u_h}(u_h(x),x)=0,
\]
for every $x\in \overline\Omega$.
\end{corollary}
\begin{proof}
The corollary follows from Proposition \ref{prop_T} and the Schauder's Fixed Point Theorem; see, for instance, \cite[Theorem 3, Section 9.2.2]{Eva1998}.
\end{proof}

Proposition \ref{prop_T} ensures $T_h$ has a fixed point in $\mathcal{B}_h$, which is tantamount to the existence of $u_h\in \mathcal{B}_h$ satisfying
\begin{equation}\label{eq_Gpf}
	G_{h}^{u_h}(u_h(x),x)=0
\end{equation}
for every $x\in \Omega_h$. Because Propositions \ref{prop_monotone2}, \ref{prop_consistency2} and \ref{prop_stability2} hold uniformly in $\varepsilon$ and $v$, we conclude the method in \eqref{eq_Gpf} is monotone, stable and consistent with \eqref{eq_reg2b} with $v$ replaced by the fixed point $u^\varepsilon$. Therefore, we are in a position to prove Theorem \ref{thm_conv2}.

\begin{proof}[Proof of Theorem \ref{thm_conv2}]
Let $G_h$ be defined in \eqref{eq_method2b}. By combining Propositions \ref{prop_monotone2}, \ref{prop_consistency2} and \ref{prop_stability2}, we conclude that $G_h$ is monotone, consistent with \eqref{eq_reg2b}, and stable. As before, define $\overline u,\underline u:\overline\Omega\to\mathbb{R}$ as
\[
	\overline u(x)\coloneqq\lim_{h\to 0}\sup_{\substack{y\to x\\y\in\Omega_h}}u_h^\varepsilon(y)\hspace{.2in}\mbox{and}\hspace{.2in}\underline u(x)\coloneqq\lim_{h\to 0}\sup_{\substack{y\to x\\y\in\Omega_h}}u_h^\varepsilon(y),
\]
where $(u_h^\varepsilon)_{0<h\ll1}$ is a family of solutions to \eqref{eq_Gpf}, for fixed $\varepsilon>0$. Arguing as in the proof of Theorem \ref{thm_convergencepure}, an application of \cite[Theorem 2.1]{BarSou} ensures that
\begin{equation}\label{eq_conv3a}
	\begin{cases}
		\varepsilon \overline u+F(D^2\overline u)-\dfrac{f}{\left(\varepsilon+|D\overline u|^2\right)^\frac{\theta^{\overline u}_\varepsilon(x)}{2}}\leq0&\hspace{.2in}\mbox{in}\hspace{.1in}\Omega\\
		\varepsilon \underline u+F(D^2\underline u)-\dfrac{f}{\left(\varepsilon+|D\underline u|^2\right)^\frac{\theta^{\underline u}_\varepsilon(x)}{2}}\geq0&\hspace{.2in}\mbox{in}\hspace{.1in}\Omega.
	\end{cases}	
\end{equation}
It remains to verify that $\overline u\leq \underline u$ on the boundary $\partial\Omega$. It follows along the same lines as in the proof of Theorem \ref{thm_convergencepure}, once one recalls that $\theta_1\leq \theta^{\overline u}_\varepsilon(x)\leq\theta_2$ and $\theta_1\leq  \theta^{\underline u}_\varepsilon(x)\leq\theta_2$ for every $x\in\overline\Omega$. This fact combines with \eqref{eq_conv3a} and builds upon the comparison principle \cite[Proposition 4]{HuaPimRamSwi} to produce $\overline u=\underline u$ in $\overline\Omega$. The convergence of the method then follows.

Finally, by letting $\varepsilon\to 0$, $u^\varepsilon$ converges locally uniformly to $u$, a viscosity solution to the free transmission problem \eqref{eq_main2}, and the proof is complete.
\end{proof}
}

\section{Numerical experiments}\label{sec_numerexp}

In this section, we present one-dimensional examples to validate the convergence of our numerical methods. We compute approximate solutions to the regularised problems \eqref{eq_reg1} and \eqref{eq_reg2}. These in turn converge to the solution of the original problems \eqref{eq_main1} and \eqref{eq_main2}, respectively, as $\varepsilon \rightarrow 0$.

Specifically, we compute an approximate solution  $u_h^\varepsilon$ to the regularised problems \eqref{eq_reg1} and \eqref{eq_reg2}, and compare it with the exact solution $u$ of the respective main problems   (\ref{eq_main1}) and  (\ref{eq_main2}). For simplicity, as in the previous sections, we omit the superscript $\varepsilon$  in the notation and refer to the numerical solutions simply as $u_h$. The following examples are constructed so that the exact solution is known a priori, allowing for a direct assessment of accuracy.

We discretise the domain using grid points $x_i=x_0+ih$, for $i=0,1,\dots, N-1,N$. Here, $x_0$ and $x_N$ denote the left and right boundaries, respectively. The finite difference approximation to the solution at each grid point $x_i$ is denoted with $u_{h,i}$. 

We recall the definitions of the discrete operators used to approximate first and second-order derivatives. The approximation for $|Du|^2$ is given by
$$
	|D_h u_{h,i}|^2 = \frac{1}{h^2}\max\left(u_{h,i}-u_{h,i-1}, u_{h,i}-u_{h,i+1},0\right)^2,
$$
which is accurate of order {\color{r2}$O(h)$} for $C^2$-regular functions. We also define the approximation of the second-order derivative as
$$
	D_h^2 u_{h,i}= \frac{1}{h^2}(u_{h,i+1}-2u_{h,i}+u_{h,i-1}).
$$
Notice the latter is accurate to order $O(h^2)$, for $C^4$-regular functions.

{\color{r2} We start by addressing the problem of the pure equation. For the pure equation the numerical scheme is 
\begin{equation*}
	G_h(u_{h,i}, x_i) : = (\varepsilon u_{h,i} +F_h(D_h^2 u_{h,i})) - \frac{f_i}{(\varepsilon + |D_h u_{h,i}|^2)^{\theta/2}}, 
\label{nm1}
\end{equation*}
for  $i=1,\dots, N-1$ and $f_i = f(x_i)$.}
 
The operator $G_h $ is also defined on the boundary as
\begin{equation*}
G_h( u_{h,i}, x_i) : = u_{h,i} -g_i,
\label{nm2}
\end{equation*}
for $i=0$ and $i=N$. Here, $g_i$ is defined in the obvious way.
This discretisation leads to a nonlinear system of equations that must be solved to obtain the approximate solution over the discrete domain. To address this, we apply the well-known explicit Euler method. It solves iteratively the nonlinear system represented by $G_h( u_{h,i}, x_i)$, for  {\color{r2}$i=1,\dots,N-1$}. See \cite{Obe}; see also \cite{BarSou}.

For $\rho>0$, define the Euler map $S_\rho(u_h)\coloneqq \left(S_\rho(u_{h,0})\ldots,S_\rho(u_{h,N})\right)$, where
\begin{equation*}
S_\rho(u_{h,i})\coloneqq u_{h,i} -\rho G_h( u_{h,i}, x_i),
\label{fpi}
\end{equation*}
for $i=0,1,\ldots,N$. The map $S_\rho$ is a contraction in the $\ell_\infty$-norm provided we choose $0<\rho\ll1$ small enough. That is, the CFL condition found before ensures a contraction; see Remark \ref{rem_cfl}. As a consequence, there exists a unique fixed point $u^*_h$, satisfying $u^*_h= S_\rho(u_h^*)$.

In the next examples, we consider $F(D^2 u)=-D^2 u$. Starting with an ansatz $u_h^{0}(x_i)$, for $i=1,\dots, N-1$, where $u_h^{0}(x_0)=u(x_0)$ and $u_{h}^{0}(x_N)=u(x_N)$ we seek the discrete solution $u_h^{n}$ through the iterative process $u_h^{n}= S_\rho u_h^{n-1}$.

We run the experiments until an (artificial) fixed terminal instant $T$ is reached. This means the number of iterations $n$ is chosen such that  $n\rho=T$. As we refine the spatial step $h$, the number of iterations required increases accordingly.
 
{\color{r1}At this point, we must also discuss the parameter $\varepsilon$ from an experimental perspective; c.f. Remark \ref{rem_he}.} Our goal is to compare the computed results with the solution $u$ of the original problem. To accurately capture the features of the solution, it is necessary that $\varepsilon \geq h$. Therefore, as $\varepsilon$ decreases, $h$ must decrease accordingly. In the following experiments, we set $\varepsilon = h$ to ensure that this condition is satisfied.
 {\color{r1} A numerical example illustrating this relationship is presented in Example~\ref{ex_pe1}.}

 We begin with two test cases for the pure problem \eqref{eq_main1}, followed by two test cases for the transmission problem \eqref{eq_main2}.

\begin{figure}[tbhp]
\center
\includegraphics[scale=0.25]{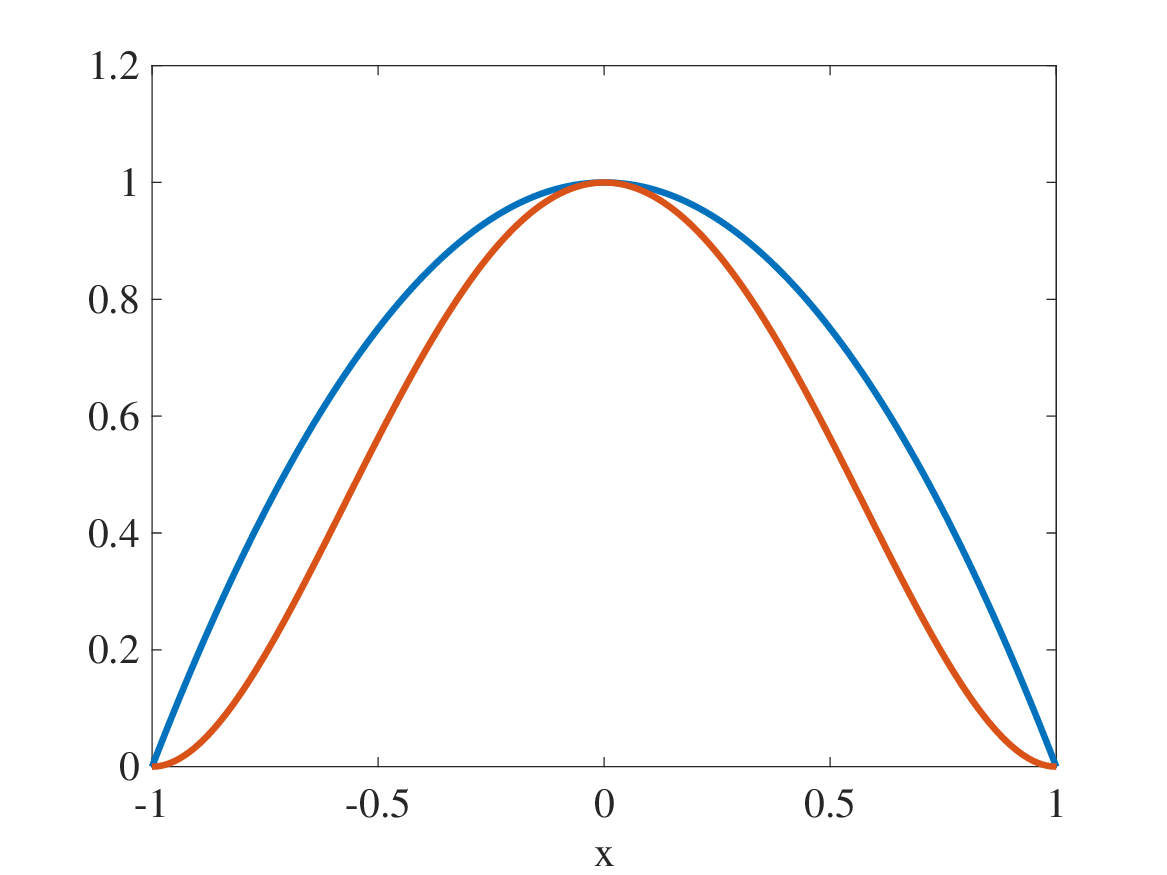}
\includegraphics[scale=0.25]{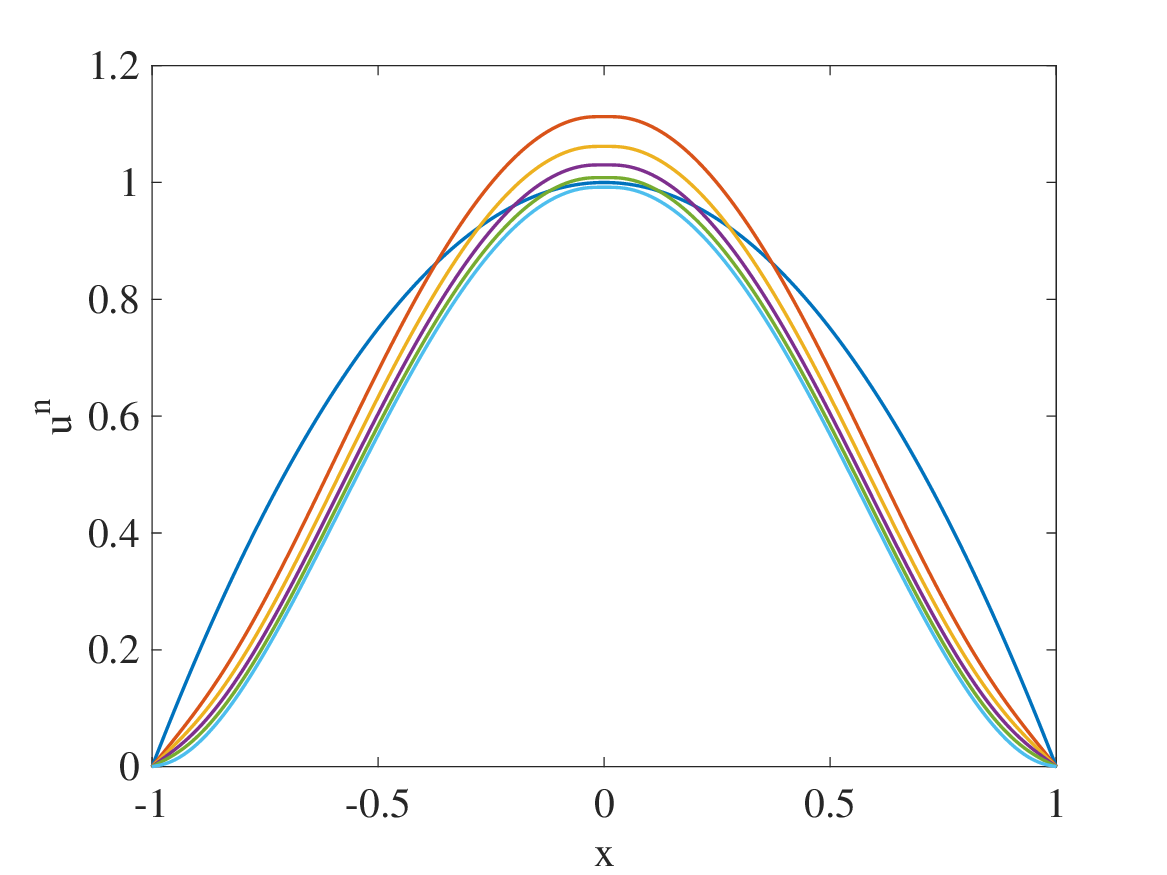}
\caption{ Example \ref{ex_pe1}: (a) Exact solution of the pure problem  (\ref{eq_main1}) (red line) and the ansatz $u_h^0$ (blue line);
(b) Solutions obtained via Euler iteration with  $h=0.01$ and $n\rho=1$. The maximum error is
$\max_i |u_{h,i} -u_i| = 8.2137\times 10^{-3}$.}
\label{F1}   
\end{figure}

\begin{example}[Pure equation I]\label{ex_pe1} \normalfont
In the first example, we consider the domain $\Omega=(-1,1)$ and the degeneracy rate $\theta=2$. The source term in problem (\ref{eq_main1}) is chosen so that the exact solution is $u(x)=(1-x)^2(1+x)^2$. The corresponding boundary conditions are $u(-1)=0$ and $u(1)=0$. In Figure \ref{F1}(a) we display the initial approximation, $u^0(x)=(1-x)(x+1)$, used to start the Euler iteration alongside the exact solution $u(x)$. Figure \ref{F1}(b) illustrates a sequence of Euler iterations performed until the approximation error falls around $h$.

The Euler iteration mimics a time-dependent process, with the parameter $\rho$ playing a role analogous to a time-step. To satisfy the stability condition, $\rho$ is chosen to depend on the mesh size $h$, specifically set as $\rho = 0.01 h^2$. For this example, we simulate until the (artificial) final time $T=1$.
\end{example}

 \begin{table}[h]
\centering
\begin{tabular}{c c | c c | c c}
\hline
\multicolumn{2}{c}{$h=0.04$} & \multicolumn{2}{c}{$h=0.01$} & \multicolumn{2}{c}{$h=0.005$} \\
\hline
$\varepsilon $ & $||\cdot||_{\infty,h}$ & $\varepsilon$ & $||\cdot||_{\infty,h}$ & $\varepsilon$ & $||\cdot||_{\infty,h}$  \\\hline
$0.01$ & $5.7619e-02$ & $0.005$ & $9.4718e-04$ & $0.001$ & $2.4970e-03$  \\
$0.02$ & $2.8996e-02$ & \cellcolor{wethepeople}$0.01$ & $8.1399e-04$ &  \cellcolor{wethepeople}$0.005$& $2.1950e-03$ \\
\cellcolor{wethepeople}$0.05$ & $2.5482e-02$ & $0.05$ & $6.1486e-03$ & $0.05$ & $4.4726e-03$ \\
$0.06$ & $2.5808e-02$  & $0.1$ & $1.1983e-02$ & $0.1$  & $1.0305e-02$ \\
$0.08$ & $2.6106e-02$  & $0.5$ & $4.3329e-02$& $0.5$  & $4.2531e-02$ \\
\hline
\end{tabular}
\caption{{\color{r1}Experimental optimal values of $\varepsilon$ with respect to $h$, where $||u_h-u||_{\infty,h}\coloneqq\max_i |u_{h,i}-u_i|$.}}
\label{table1}
\end{table}

{\color{r1}
We also use this example to look for experimental information on the connection between the regularisation and the discretisation parameters. In Table~\ref{table1}, we present the experimental optimal values of $\varepsilon$ in relation to $h$.
The discrete error displayed is defined by $||u_h-u||_{\infty,h}\coloneqq\max_i |u_{h,i}-u_i|$, where $u_h$ is the approximate solution of problem \eqref{eq_reg1} and $u$ is the exact solution of problem \eqref{eq_main1}.
 Our findings suggest that as $\varepsilon$ decreases, so does the error. However, as $\varepsilon$ crosses the $h$-value threshold, the approximation ceases to improve, and may even worsen. Consequently this indicates that  smaller values of $\varepsilon$ require mesh refinement to achieve error reduction efficiently.}
 
\begin{figure}[tbhp]
\center
\includegraphics[scale=0.25]{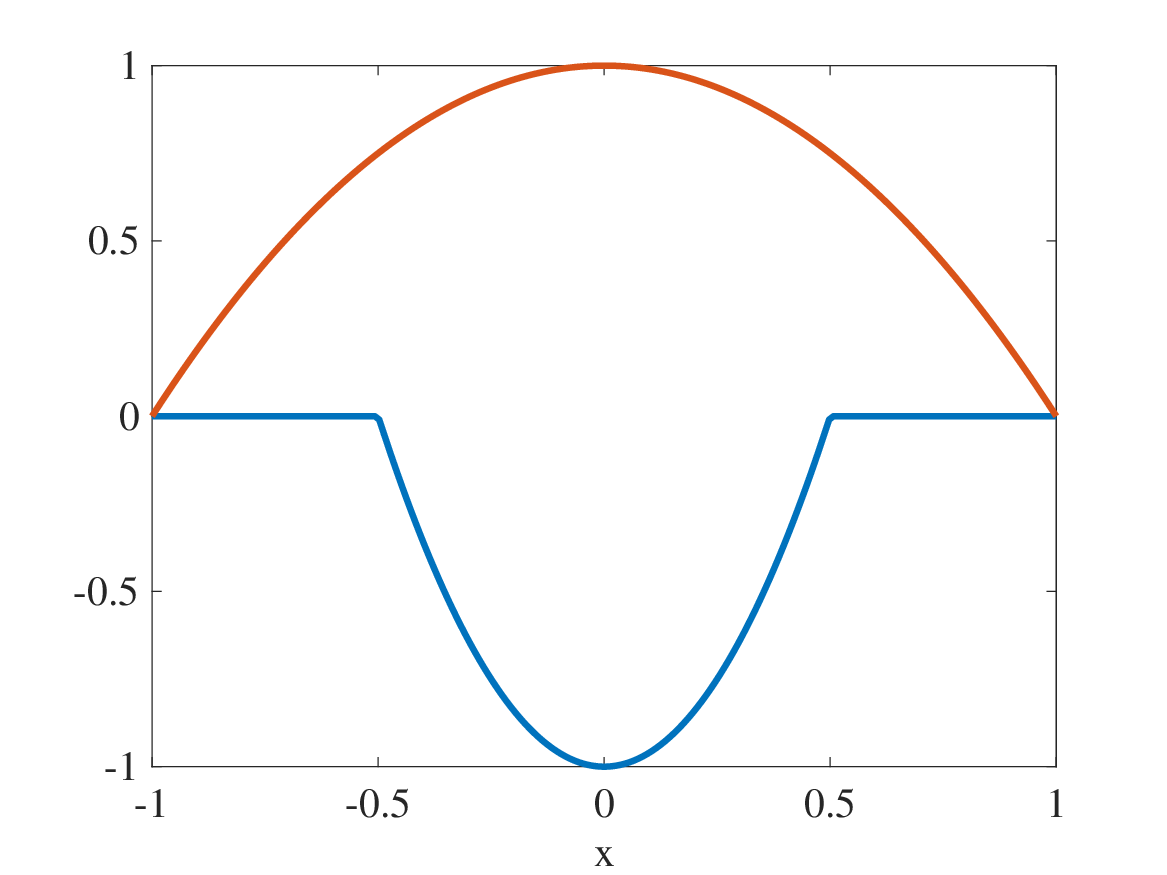}
\includegraphics[scale=0.25]{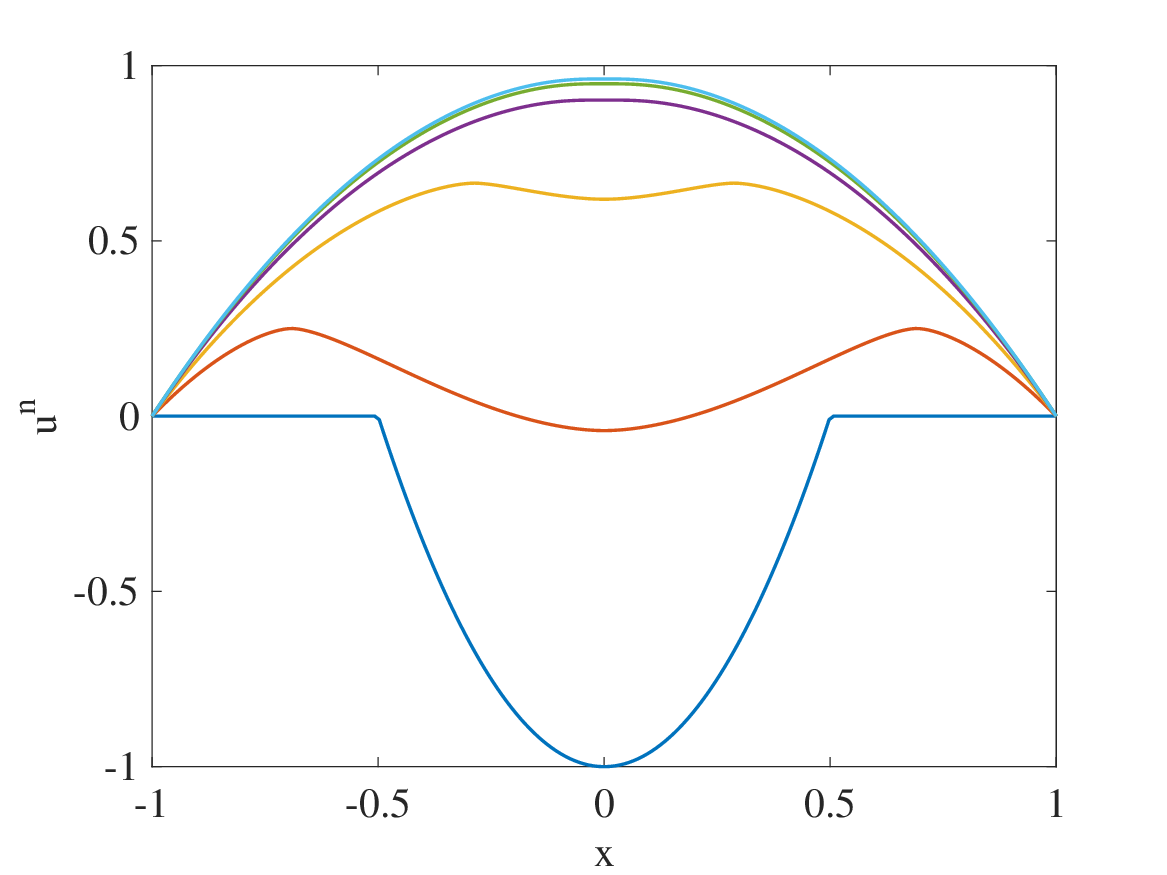}
\caption{ Example \ref{ex_pe2}: (a) Exact solution of the pure problem  (\ref{eq_main1}) (red line) and the ansatz $u_h^0$ (blue line);
(b) Solutions obtained via Euler iteration with  $h=0.01$ and $n\rho=1$. The maximum error is
$\max_i |u_{h,i} -u_i| = 3.7964\times 10^{-2}$.}
\label{F2}   
\end{figure}

\begin{example}[Pure equation II]\label{ex_pe2} \normalfont
Our second example is set in the interval $\Omega=(-1,1)$, with a degeneracy rate $\theta=1$. The exact solution in this case is $u(x)=1-x^2$, which forces $u(-1)=u(1)=0$. We start the Euler method with an ansatz whose gradient vanishes in the set $(-1,-1/2)\cup(1/2,1)$. A relevant aspect of the present example concerns the regions $\{|Du|=0\}$. Heuristically, the equation provides no information on $u$ along this region. However, the regularisation in the method is capable of capturing the structure of the PDE even in this case. Figure \ref{F2}(a) showcases the ansatz and the exact solution $u(x)$. Figure \ref{F2}(b) presents a sequence of Euler iterations yielding an approximation error of order $h$. Once again,  $\rho = 0.01 h^2$ and $T=1$.
\end{example}

When it comes to the free transmission problem, we present two examples. In the first one, the exact solution has a strictly positive derivative, whereas in the second one, the gradient of the solution vanishes at $x=0$. We notice our method is capable of bypassing the effects of the gradient-driven degeneracy in both cases.

{\color{r2} 
For the transmission problem, the exponent $\theta$ is replaced by the function $\theta_\varepsilon^v(x)$. 
Consequently, the numerical method requires an additional iteration to compute the fixed point associated with $v$. The resulting numerical scheme is given by
\begin{equation*}
	G_h^{v_{n}}(u_{h,i}, x_i) : = (\varepsilon u_{h,i} +F_h(D_h^2 u_{h,i})) - \frac{f_i}{(\varepsilon + |D_h u_{h,i}|^2)^{\theta_i^{v_n}/2}}, 
\label{nm1}
\end{equation*}
for  $i=1,\dots, N-1$, where $\theta_i^{v_n} \coloneqq \theta^{v_n}(x_{i})$ and $f_i = f(x_i)$.

More specifically, for a fixed $\varepsilon>0$, we begin with an initial guess $v_1$. At the $n$th outer iteration, given the current approximation $v_n$, we solve the nonlinear discrete problem
\begin{equation}
G_h^{v_n}(u_h^n,x_i)=0,\qquad x_i\in\Omega_h,
\label{gknp}
\end{equation}
to obtain the discrete solution $u_h^n$. We then update the iterate by setting
\[
v_{n+1}=u_h^n,
\]
and repeat the procedure using the updated function $v_{n+1}$. The iteration is continued until $n=n_{\max}$.
The final iterate $u_h^{n_{\max}}$ is taken as the numerical approximation of the solution of
\[
G_h^{u_h}(u_h,x)=0.
\]
For each outer iteration $n$, equation~\eqref{gknp} gives rise to a nonlinear system whose unknowns are the values of $u_h^n$ at the interior grid points.
%
This operator also depends on   boundary points. However,  since their values are known, they are not treated as unknowns in the system.
To obtain the approximate solution over the discrete domain, the resulting nonlinear system of equations must be solved. 
For each fixed $v_n$, we apply the Euler map. Starting from an initial guess $u_h^{n,0}$, we perform $m_{\max}$ iterations according to
$$
u_{h,i}^{n,m+1}  =u_{h,i}^{n,m} - \rho G_{h}^{v_{n}}(u_{h,i}^{n,m},x_i), \quad m=0,1,\dots, m_{max}-1.
$$
Therefore for each $v_n$ we obtain the numerical solution $u_{h}^{n,m_{max}}$ and at the end we obtain
the numerical solution $u_{h}^{n_{max},m_{max}}$.
  }

\begin{figure}[h]
\center
\includegraphics[scale=0.17]{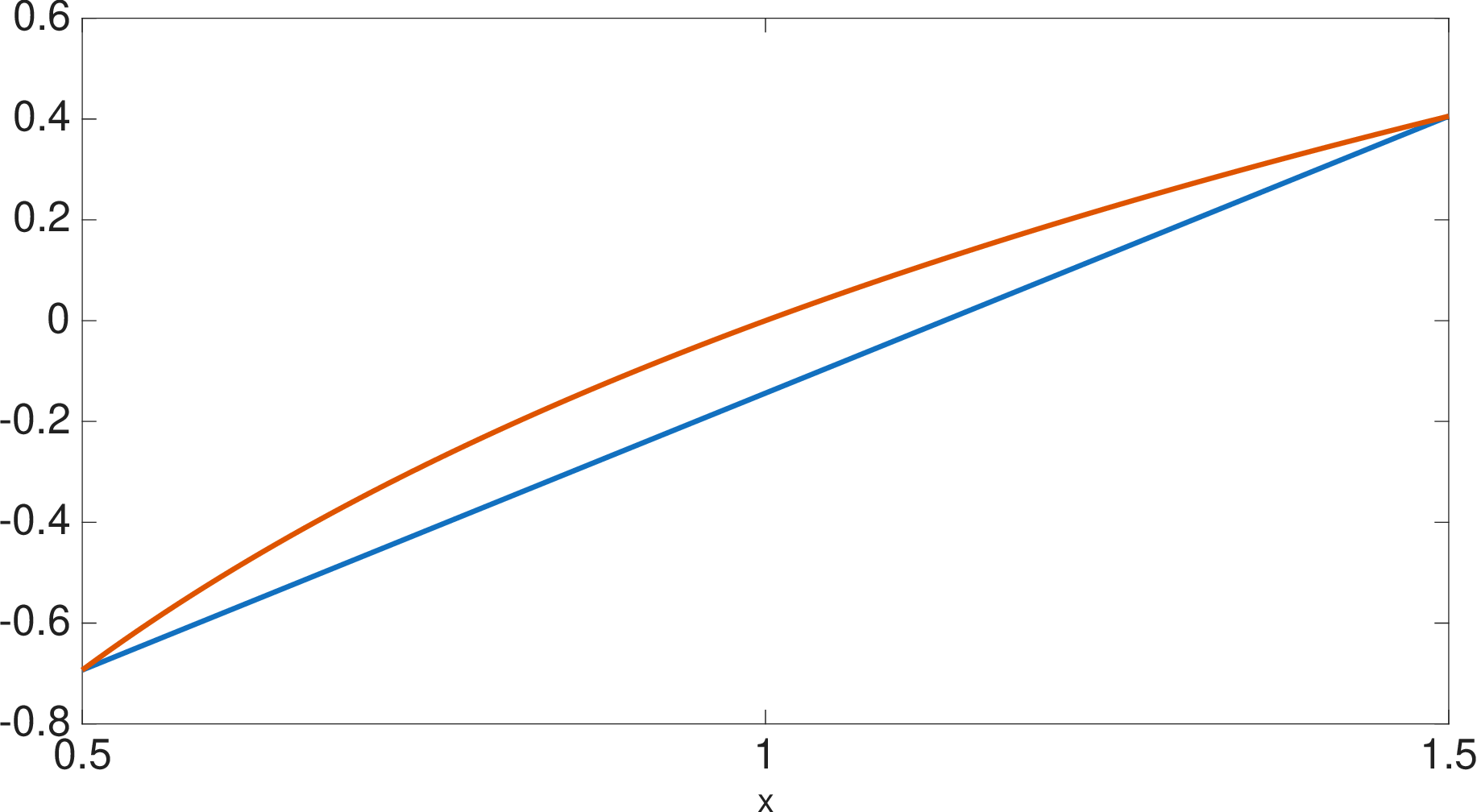}
\hspace*{1cm}
\includegraphics[scale=0.17]{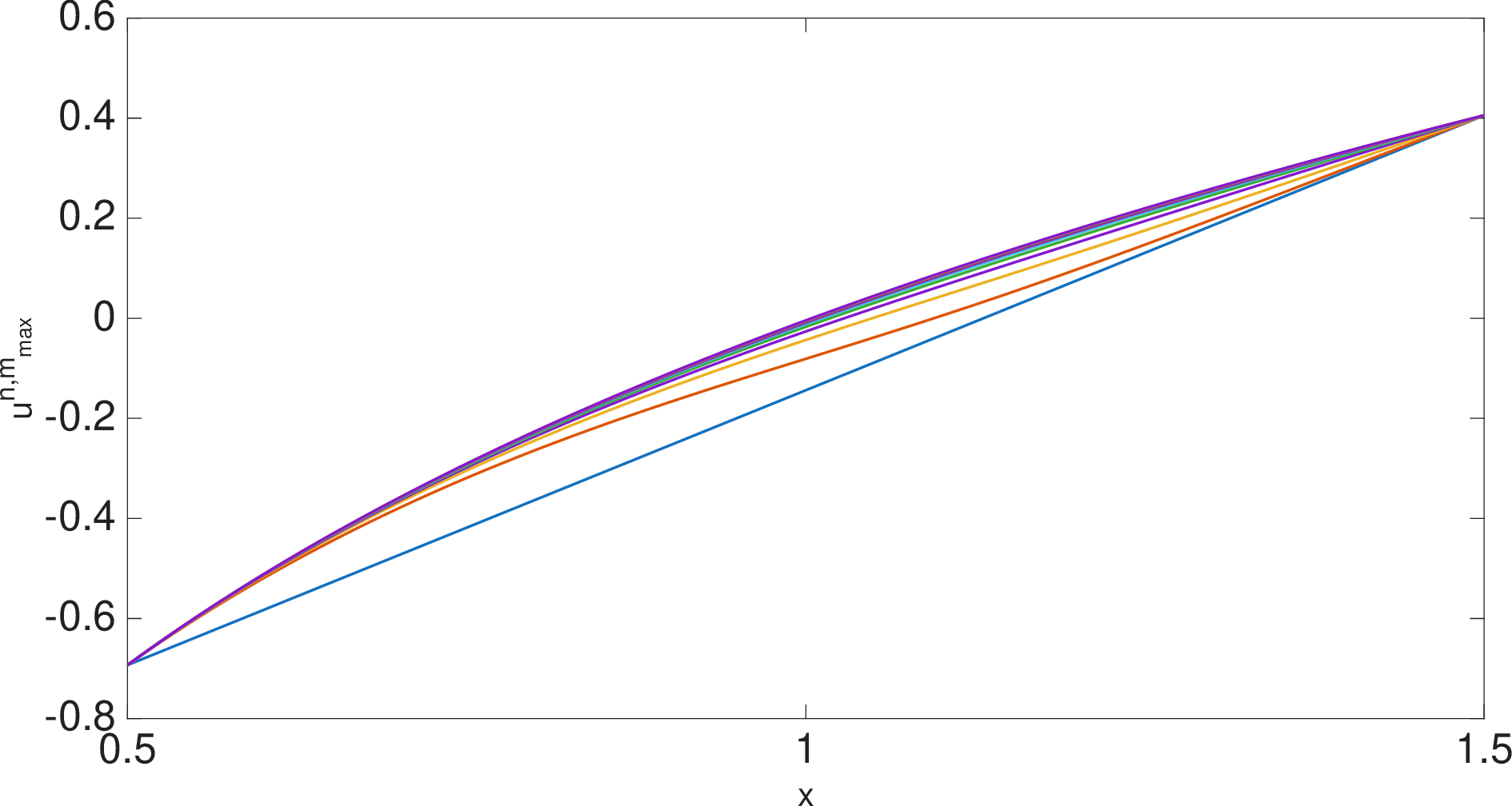}

\caption{ Example \ref{ex_tp1}: (a) Exact solution of problem  (\ref{eq_main2})  (red line) and the ansatz $u_h^{1,0}$ (blue line);
(b) Solution $u_h^{n,m_{max}}$ for some values of $n$ obtained with $h=0.01$ and $n_{max}m_{max}\rho=0.3$. 
 The maximum error is $\max_i |u_{h,i}^{n_{max},m_{max}} -u_i|  = 4.2872\times 10^{-3}$.}
\label{F3}   
\end{figure}

\begin{example}[Transmission problem I -- strictly positive gradient]\label{ex_tp1}\normalfont
Consider the domain $\Omega=(1/2,3/2)$ and the degeneracy rates $\theta_1=2$ and  $\theta_2=4$. The source term in problem (\ref{eq_main2}) is chosen so that the exact solution is $u(x)=\log(x)$. The corresponding boundary conditions are $u(1/2)=\log(1/2)$ and $u(3/2)=\log(3/2)$.
  
The initial estimate for the Euler iteration and the initial $v_1$ are
$$
v_1(x_i)=u_{h}^{1,0}(x_i)=\log(3)x_i+\log\left(\frac{1}{2\sqrt{3}}\right), \ i=1,\dots, N-1.
$$
In Figure \ref{F3}(a), we display the initial approximation, $u^{1,0}(x)$, used to start the Euler iteration alongside the exact solution $u(x)$. In Figure \ref{F3} (b), we display how Euler's method progresses until the solution approximates the exact solution with an error smaller than $h$. To satisfy the stability condition, $\rho$ is chosen to depend on the mesh size $h$, specifically set as $\rho = 0.01 h^2$. For this example, we simulate until the final time $T=0.3$. We need fewer iterations than in the previous example, due to the geometry of the exact solution.
\end{example}

\begin{figure}[h]
\center
\includegraphics[scale=0.17]{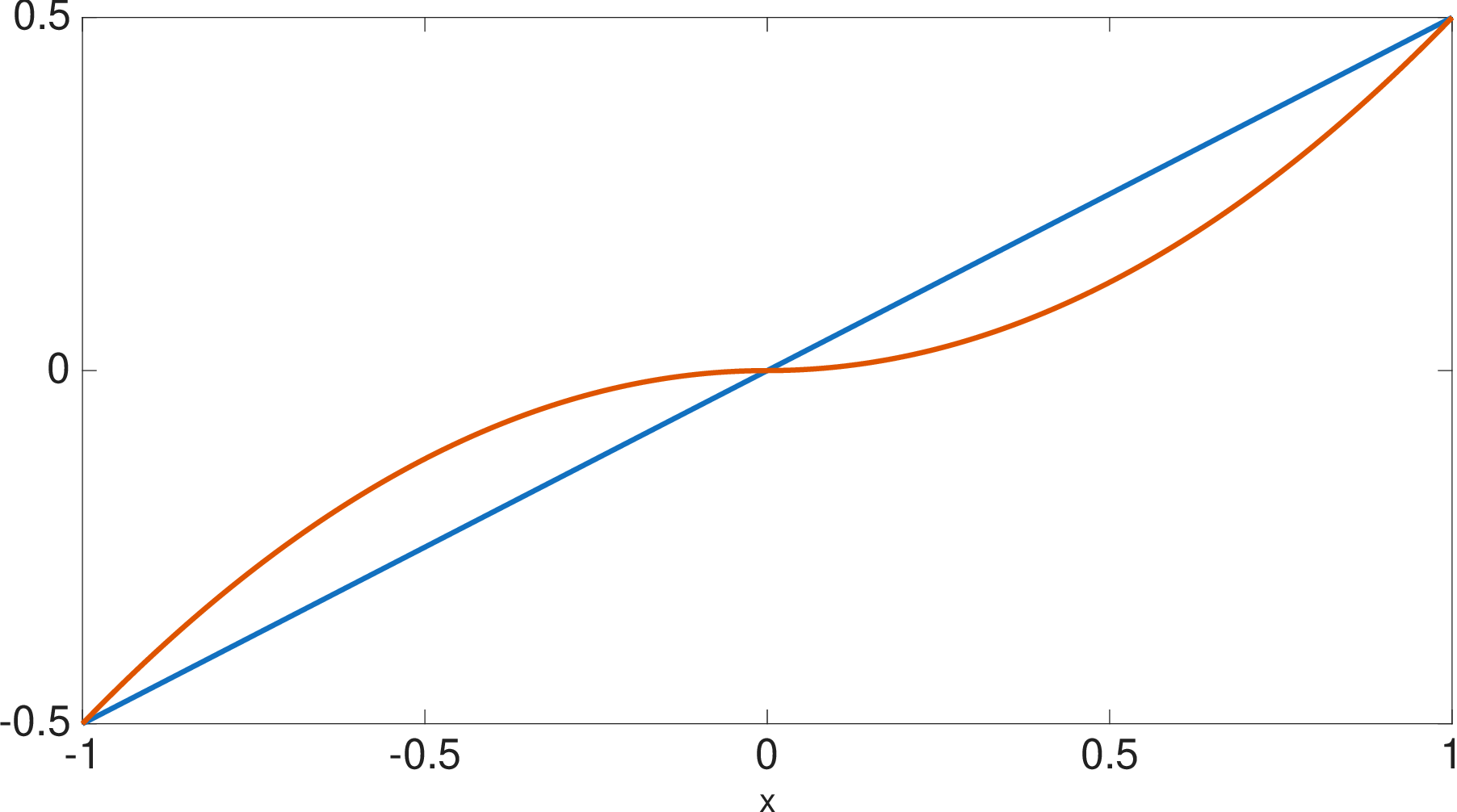}\hspace*{1cm}
\includegraphics[scale=0.17]{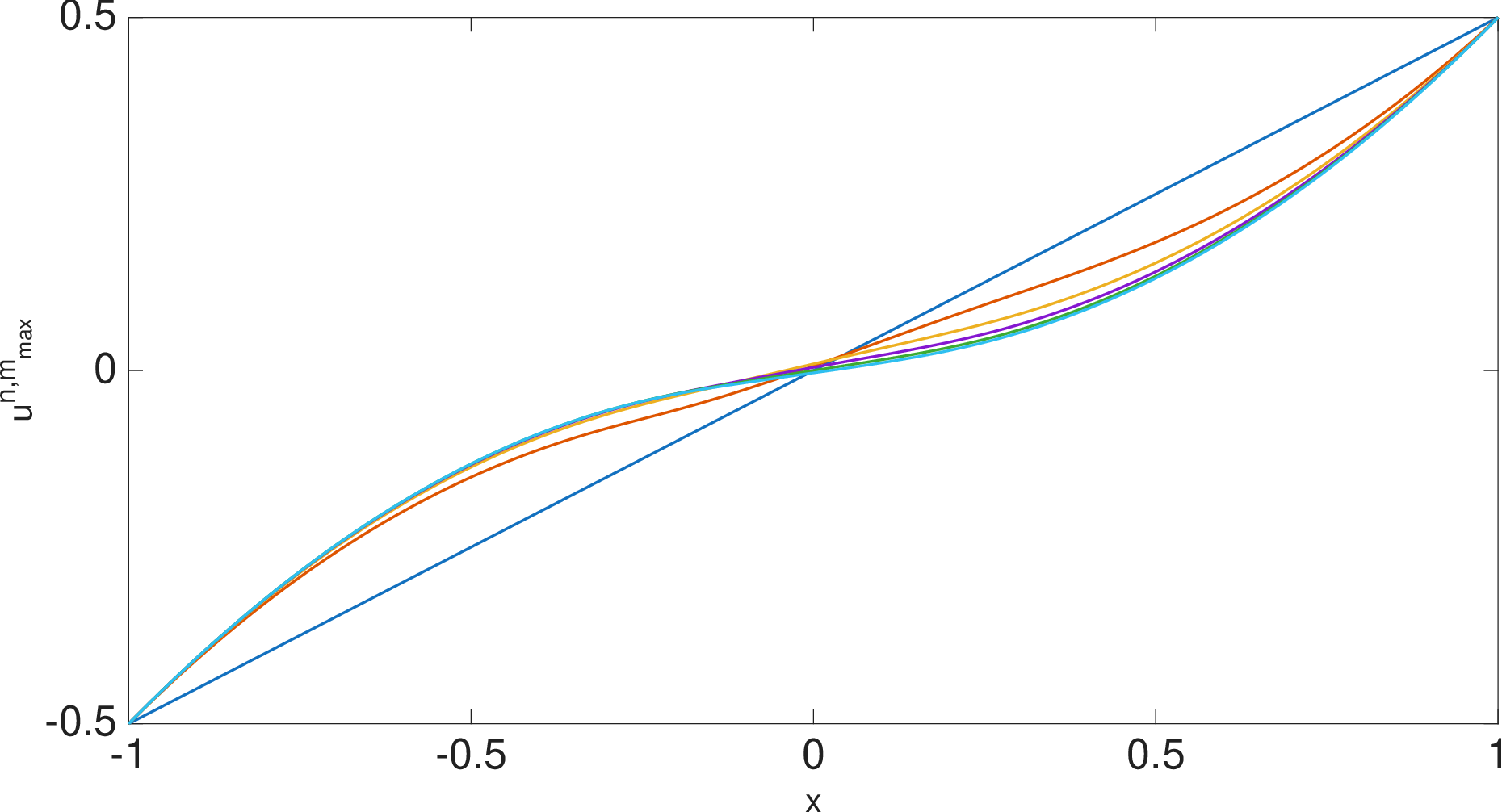}
\caption{ Example \ref{ex_tp2}:  (a) Exact solution of problem (\ref{eq_main2}) (red line) and the ansatz $u_h^{1,0}$ (blue line);
 (b) Solution $u_h^{n,m_{max}}$ for some values of $n$ obtained with $h=0.01$ and $n_{max}m_{max}\rho=0.3$. 
 The maximum error is $\max_i |u_{h,i}^{n_{max},m_{max}} -u_i| = 1.3393\times 10^{-2}$.}
\label{F4}   
\end{figure}

\begin{example}[Transmission problem II -- vanishing gradient]\label{ex_tp2} \normalfont
Here, $\Omega=(-1,1)$, $\theta_1=2$, and $\theta_2=4$. The source term in problem \eqref{eq_main2} is chosen so that the exact solution is $u(x)=x^2/2$ for $x>0$ and $u(x)=-x^2/2$ for $x<0$. The corresponding boundary conditions are $u(-1)=-1/2$ and $u(1)=1/2$. The initial estimate for the Euler iteration is 
$$
	v_1(x_i)=u_{h}^{1,0}(x_i)=\frac{x_i}{2}, \ i=1,\dots, N-1.
$$
We choose $\rho$ as in the previous examples, $\rho = 0.01 h^2$. In Figure \ref{F4}(a), we display the initial approximation, $u^{1,0}(x)$, used to start the Euler iteration alongside the exact solution $u(x)$. Figure \ref{F4}(b) shows the Euler method iterations until the solution approximates the exact solution with an error around $h$.

\end{example}

We kept the final value $n_{max}m_{max}\rho=0.3$ the same as in the previous example for comparison purposes. This fourth example requires more iterations to reach a satisfactory approximation of the solution compared to the previous one. The most challenging region to approximate appears to be around $x=0$. The maximum error is slightly larger than that in Example \ref{ex_tp1}.

\bibliographystyle{plain}
\bibliography{biblio}

\end{document}